\documentclass[titlepage,11pt,leqno]{report}
\usepackage{latexsym}
\usepackage{doublespace}
\usepackage{amsmath}
\usepackage{amsfonts}
\usepackage{threeparttable}
\usepackage[mathscr]{eucal}
\usepackage[dvips]{graphicx}
\begin{document}
\baselineskip=17.7pt
\bibliographystyle{unsrt}
\newtheorem{theorem}{\indent \sc Theorem}
\newtheorem{corollary}{\indent \sc Corollary}
\newtheorem{lemma}{\indent \sc Lemma}
\newtheorem{proposition}{\indent \sc Proposition}
\newtheorem{remark}{\indent \sc Remark}

\newcommand{\re}[0]{\mathbb{R}}
\newcommand{\ind}[0]{\mbox{\Large\bf 1}}
\newcommand{\pr}[0]{\mbox{P}}

\begin{center}
{\bf\large MARGINAL ASYMPTOTICS FOR THE ``LARGE P, SMALL N'' PARADIGM: WITH
APPLICATIONS TO MICROARRAY DATA}
\vspace{0.25in}

{\bf University of Wisconsin, Madison, Department of Biostatistics and 
Medical Informatics Technical Report 188}
\vspace{0.25in}

{\sc By Michael R. Kosorok} and {\sc Shuangge Ma}

{\em University of Wisconsin-Madison and University of Washington}
\end{center}
\vspace{0.5in}

\begin{quote}
The ``large p, small n'' paradigm arises in microarray studies,
where expression levels of thousands of genes are monitored
for a small number of subjects. There has been an increasing 
demand for study of asymptotics
for the various statistical models and methodologies using genomic data.  
In this article, we focus on one-sample and two-sample microarray
experiments, where the goal is to identify significantly 
differentially expressed genes. We establish uniform
consistency of certain estimators of marginal distribution functions,
sample means and sample medians under the large p small n 
assumption. We also establish uniform consistency of
marginal $p$-values based on certain asymptotic approximations
which permit inference based on false discovery rate techniques.
The affects of the normalization process on these results is also 
investigated. Simulation studies and data analyses are used to
assess finite sample performance.
\vspace{0.25in}

{\it Key words and phrases.} Brownian bridge, Brownian motion, Empirical 
Process, False discovery rate, Hungarian construction, Marginal 
asymptotics, Maximal inequalities, Microarrays.
\end{quote}

\newpage
{\bf 1. Introduction.}
Microarrays are capable of monitoring gene expression
on a large scale and are becoming a routine tool in
biomedical research. Studies of associations between
microarray measurements and variations of phenotypes 
can lead to a better treatment assignment and so there
has been an increasing demand for novel statistical tools
analyzing such data. For example, several recent developments in
microarray data analysis have involved semiparametric model methodology.
Such research includes, but is not limited to, 
estimation of normalization effects with a semi-linear in-slide model (SLIM) 
in Fan, Peng and Huang (2004) (FPH hereafter), estimation and inference 
of gene effects in Yang et al. (2001) and Huang, Wang and Zhang (2005)
(HWZ hereafter), classification of phenotypes based on Affymetric genechip 
data in Ghosh and Chinnaiyan (2004), and survival analysis with right censored
data and genomic covariates (Gui and Li, 2004). 

Although statistical analysis with microarray data has been one
of the most investigated areas, theoretical studies of asymptotic
properties of different statistical methodologies remain rare (for
important exceptions to this, see van der Laan and Bryan, 2001; FPH; 
and HWZ). The paucity of such research is
partly caused by the abnormal type of asymptotics associated
with microarrays: the dimension of the covariate $p$ is usually much
larger than the sample size $n$, i.e., the "large p, small n"
paradigm referred to in West (2003).
In this article, we focus on asymptotics for the simple settings
of one-sample and two-sample comparisons, where the goal is to find
genes differentially expressed for different phenotype groups.

Consider, for example, 
a simple one-sample cDNA microarray study, where the goal is to 
identify genes differentially expressed from zero. Note that this data 
setting and the following discussions can be easily extended to incorporate
two-sample microarray studies as in Yang et al. (2001).
Studies using Affymetrix genechip data can be included in the
same framework with only minor modifications.
Denote $Y_{ij}$ and $Z_{ij}$ as the background-corrected
log-ratios and log-intensities (as in HWZ), for 
array $i=1,\ldots,n$ and gene $j=1,\ldots,p$. We consider the following
simplified partial linear model for cDNA microarray data:
\begin{eqnarray}\label{model.1}
Y_{ij}=\mu_j+h_i(Z_{ij})+\epsilon_{ij},
\end{eqnarray}
where $\mu_j$ are the fixed gene effects, $h_i(Z_{ij})$ are the smooth 
array-specific normalization effects (constrained to have mean zero
within each array) and $\epsilon_{ij}$ are mean 
zero (within array) random errors. The constraints are for model
identifiability. For simplicity of exposition, we have omitted other 
potentially important terms in our model, such
as possible print-tip effects, and array-specific
position and scale constants.  We note, however, that the theory
we present in this paper can extend readily to these richer models.

Models similar to~\ref{model.1} have been investigated by HWZ
and FPH. In HWZ, asymptotic properties 
based on least squares estimation are established assuming fixed $p$ and 
$n\to \infty$. It is shown that $\mu_j$ and $h_i$ can both be consistently 
estimated with optimal convergence rates. 
In FPH, partial consistency type asymptotics are 
established. It is proved that when $n$ is fixed and $p\to \infty$, $h_i$ 
can be consistently estimated by an estimator $\hat{h}_i$, 
although $\mu_j$ cannot be consistently estimated. 
If we let $X_{ij}=\mu_i+\epsilon_{ij}$ and $\tilde{X}_{ij}
=Y_{ij}-\hat{h}_i(X_{ij})$, the results of FPH can be restated
as $\max_{1\leq i\leq n}\max_{1\leq j\leq p}|\tilde{X}_{ij}-X_{ij}|=o_P(1)$.
In otherwords, the normalization process is consistent. This 
permits the use of the normalized array-specific
gene effects $\tilde{X}_{ij}$ for inference in place of the
true array-specific gene effects $X_{ij}$. However, because $n$ is
fixed, the permissible inference tools at the gene level are restricted 
to exact methods, such as permutation tests.

The goal of our paper is to study normalization and inference when
the number of arrays $n\rightarrow\infty$ slowly while the number of
genes $p>>n$. This is essentially the same asymptotic framework
considered in van der Laan and Bryan (2001) who show that provided
the range of expression levels is bounded, the sample means 
consistently estimate the mean gene effects uniformly across genes 
whenever $\log p=o(n)$. We extend the results of van der Laan
and Bryan (2001), FPH and HWZ in three important ways. First, uniform 
consistency results are extended to general empirical distribution 
functions and sample medians. Second, a precise Brownian bridge
approximation to the empirical distribution function is developed
and utilized to establish uniform validity of marginal p-values based 
on approximations which are asymptotic in $n$. The statistical tests we 
consider for this purpose include both one and two sample mean and
median tests as well as several other functionals of the empirical 
distribution function. We find that the rate requirement is either
$\log p_n=o(n^{1/2})$ or $\log p_n=o(n^{1/3})$, depending on the
choice of test statistic. Third, these results are 
further extended to allow for the presence of normalization error.

An important consequence of these results is that approximate
p-values based on normalized gene expression data can be validly
applied to false discovery rate (FDR) methods
(see Benjamini and Hochberg, 1995) for identifying differentially
expressed genes. We refer to this kind of asymptotic regime as
``marginal asymptotics'' (see also Kosorok and Ma, 2005) because
the focus of the inference is at the marginal (gene) level,
even though the results are uniformly valid over all
genes. The main conclusion of our paper is that the marginal asymptotic
regime is valid even if the number of genes increases
almost exponentially relative to the number of arrays, i.e.,
$\log p_n=o(n^{\alpha})$ for some $\alpha>0$. Qualitatively,
this seems to be the correct order of asymptotics for microarray
experiments with a moderate number, say $\sim 50$, of replications.
The main tools we use to obtain these results include maximal
inequalities, a specialized Hungarian construction for the
empirical distribution function, and a precise bound on the
modulus of continuity of Brownian motion.

The article is organized as follows. In sections 2--4, we investigate
marginal asymptotics based on the true gene effects (no normalization
error). Section~2 discusses one-sample inference based on the mean
and the median. Section~3 extends section~2 to the two-sample setting.
Section~4 considers one and two sample inference when the statistics
are distribution free. Section~5 demonstrates under reasonable regularity 
conditions that the asymptotic results obtained in sections~2--4 are not 
affected by the normalization process. Simulation studies and data analyses 
in section 6 are used to assess the finite sample performance 
and to demonstrate the practical utility of the proposed asymptotic theory.
A brief discussion is given in section 7. Proofs are given in section 8.   

\vspace{0.35cm}
{\bf 2. Marginal asymptotics for one sample studies}. 
The results of this section are based on the true data (without
normalization error). For each $n\geq 1$, 
let $X_{1(n)},\ldots,$ $X_{n(n)}$ 
be a sample of i.i.d. vectors of length $p_n$, 
where the dependence within vectors
is allowed to be arbitrary. Denote the $j$th component of the $i$th
vector $X_{ij(n)}$, i.e., $X_{i(n)}=(X_{i1(n)},\ldots,X_{ip_n(n)})'$.
Also let the marginal distribution of $X_{1j(n)}$ be denoted
$F_{j(n)}$, and let $\hat{F}_{j(n)}(t)=n^{-1}\sum_{i=1}^n
\ind\{X_{ij(n)}\leq t\}$, for all $t\in\re$ and each $j=1,\ldots,p_n$,
where $\ind\{A\}$ is the indicator of $A$. Note that $n$ can be viewed
as the number of microarrays while $p_n$ can be viewed as the number 
of genes. As mentioned in the introduction, our asymptotic interest
focuses on what happens when $n$ increases slowly while $p_n$ 
increases rapidly.

We first establish, in section~2.1, uniform consistency of the marginal 
empirical distribution function estimator and also the uniformity of a 
Brownian bridge approximation to the standardized version of this estimator.
These results are then used in section~2.2 to establish uniform
consistency of the marginal sample means and uniform validity of marginal
p-values based on the normal approximation to the t-test. The results
are extended in section~2.3 for inference based on the marginal sample 
medians. Note that both the mean and median are functionals of the
empirical distribution function. The mean is computationally simpler,
but the median is more robust to data contamination.
 
\vspace{0.35cm}
{\it 2.1 Consistency of the marginal empirical distribution functions}.
The results of this section
will form the basis for the results presented in sections~2.2 and~2.3.
The two theorems of this section, theorems~\ref{t.c1} and~\ref{t1} below,
are somewhat surprising, high dimensional extensions of two classical 
univariate results 
for empirical distribution functions: the celebrated Dvoretsky, Kiefer
and Wolfowitz (1956) inequality as refined by Massart (1990) and the 
celebrated Koml\'{o}s, Major and Tusn\'{a}dy (1976) Hungarian construction 
as refined by Bretagnolle and Massart (1989). The extensions utilize maximal
inequalities based on Orlicz norms (see chapter~2.2 of van der Vaart
and Wellner, 1996). For any real random variable $Y$ and any $d\geq 1$, 
let $\|Y\|_{\psi_d}$ denote the Orlicz norm for $\psi_d(x)=e^{x^d}-1$, i.e.,
$\|Y\|_{\psi_d}=\inf\left\{C>0:\mbox{E}\left[e^{|Y|^d/C}-1\right]\leq 1
\right\}$. Note that these norms increase with $d$ (up to a constant
depending only on $d$) and that $\|\cdot\|_{\psi_1}$ dominates
all $L_p$ norms (up to a constant depending only on $p$). Also let
$\|\cdot\|_{\infty}$ be the uniform norm.

The first theorem we present yields simultaneous consistency
of all the $\hat{F}_{j(n)}$s for the corresponding $F_{j(n)}$s:
\begin{theorem}\label{t.c1}
There exists a universal constant $0<c_0<\infty$ such that, 
for all $n,p_n\geq 2$,
\begin{eqnarray}
\left\|\max_{1\leq j\leq p_n}\left\|\hat{F}_{j(n)}-F_{j(n)}
\right\|_{\infty}\;\right\|_{\psi_2}&\leq& 
c_0\sqrt{\frac{\log p_n}{n}}.\label{c1.e1}
\end{eqnarray}
In particular, if $n\rightarrow\infty$
and $\log p_n/n\rightarrow 0$, then the left-hand-side of~(\ref{c1.e1})
goes to zero.
\end{theorem}

\begin{remark}\label{r.t.c1}
One can show that the rate on the right-side of~(\ref{c1.e1}) is
sharp, in the sense that there exist sequences of
data sets, where $(\log p_n/n)^{-1/2}\max_{1\leq j\leq p_n}
\|\hat{F}_{j(n)}-F_{j(n)}\|_{\infty}\rightarrow c$, in probability,
as $n\rightarrow\infty$, and where $0<c<\infty$. In particular, the
statement is true if the genes are all independent, 
$n,p_n\rightarrow\infty$ with $\log p_n=o(n)$, and $c=1/2$.
\end{remark}

The second theorem shows that the standardized empirical
processes $\sqrt{n}(\hat{F}_{j(n)}-F_{j(n)})$ can be simultaneously
approximated by Brownian bridges in a manner which preserves the
original dependency structure in the data. This feature will be useful 
in studying FDR (see Benjamini and Hochberg, 1995) properties later on.
To this end, let ${\cal F}_{j(n)}$ 
denote the smallest $\sigma$-field making all
of $X_{1j(n)},\ldots,X_{nj(n)}$ measurable, $1\leq j\leq p_n$. Also
let ${\cal F}_n$ be the smallest $\sigma$-field making all
of ${\cal F}_{1(n)},\ldots,{\cal F}_{p_n(n)}$ measurable. 
\begin{theorem}\label{t1}
There exists universal constants $0<c_1,c_2<\infty$ such that, for
all $n,p_n\geq 2$,
\begin{eqnarray}
\left\|\max_{1\leq j\leq p_n}\left\|\sqrt{n}(\hat{F}_{j(n)}-F_{j(n)})-
B_{j(n)}(F_{j(n)})\right\|_{\infty}\;\right\|_{\psi_1}&\leq& \frac{c_1
\log n+c_2\log p_n}{\sqrt{n}},\label{t1.e1}
\end{eqnarray}
for some stochastic processes
$B_{1(n)},\ldots,B_{p_n(n)}$ which are conditionally independent given
${\cal F}_n$ and for which each
$B_{j(n)}$ is a standard Brownian bridge with conditional
distribution given ${\cal F}_n$ depending only on ${\cal F}_{j(n)}$,
$1\leq j\leq p_n$.
\end{theorem}

\vspace{0.35cm}
{\it 2.2 Estimation of marginal sample means}.
Now we consider marginal inference based on the marginal
sample mean. For each $1\leq j\leq p_n$, assume for this section
that the closure
of the support of $F_{j(n)}$ is a compact interval $[a_{j(n)},b_{j(n)}]$
with $a_{j(n)}\neq b_{j(n)}$, and that $F_{j(n)}$ has
mean $\mu_{j(n)}$ and standard deviation $\sigma_{j(n)}>0$. 
Let $\bar{X}_{j(n)}$ be the sample mean of $X_{1j(n)},\ldots,X_{nj(n)}$.
The following corollary yields simultaneous consistency of the marginal
sample means:  
\begin{corollary}\label{c3}
Under the conditions of theorem~\ref{t.c1} and with the same constant
$c_0$, we have for all $n,p_n\geq 2$,
\begin{eqnarray}
\left\|\max_{1\leq j\leq\infty}|\bar{X}_{j(n)}-\mu_{j(n)}|\;\right\|_{\psi_2}
\leq c_0\sqrt{\frac{\log p_n}{n}}
\max_{1\leq j\leq p_n}|b_{j(n)}-a_{j(n)}|.\label{c3.e1}
\end{eqnarray}
\end{corollary}

\begin{remark}\label{r3}
Note that corollary~\ref{c3} slightly 
extends the large $p$ small $n$ consistency
results of van der Laan and Bryan (2001) by allowing the range of the 
support to increase with $n$ provided it does not increase too rapidly.
\end{remark}

Now assume that we wish to test the marginal
null hypothesis $H_0^{j(n)}:\mu_{j(n)}=\mu_{0,j(n)}$ with the test statistic
\[T_{j(n)}=\frac{\sqrt{n}(\bar{X}_{j(n)}-\mu_{0,j(n)})}
{\hat{\sigma}_{j(n)}},\]
where $\hat{\sigma}_{j(n)}$ is a location-invariant and consistent  
estimator of $\sigma_{j(n)}$.
To use FDR, we need to obtain uniformly consistent estimates of the
p-values of these tests.  One way to do this is with permutation
methods. A computationally easier way is to just use
$\hat{\pi}_{j(n)}=2\Phi(-|T_{j(n)}|)$, where $\Phi$ is the distribution
function for the standard normal. The conclusion of the following
corollary is that this approach leads to uniformly consistent p-values
under reasonable conditions:

\begin{corollary}\label{c2}
Let the constants $c_1,c_2$ be as in theorem~\ref{t1}.
Then, for all $n,p_n\geq 2$, there exist standard normal random variables
$Z_{1(n)},\ldots,Z_{p_n(n)}$ which are conditionally independent
given ${\cal F}_n$ and for which each $Z_{j(n)}$ has conditional
distribution given ${\cal F}_n$ depending only on ${\cal F}_{j(n)}$,
$1\leq j\leq p_n$, such that
\begin{eqnarray}
\max_{1\leq j\leq p_n}\left|\hat{\pi}_{j(n)}-\pi_{j(n)}\right|&\leq&
\frac{c_1\log n+c_2\log p_n}{\sqrt{n}}\left(\max_{1\leq j\leq p_n}
\frac{|b_{j(n)}-a_{j(n)}|}{\sigma_{j(n)}}\right)\label{c2.e1}\\  
&&+\frac{1}{2}\left(\max_{1\leq j\leq n}
(\hat{\sigma}_{j(n)}\vee\sigma_{j(n)})
\left|\frac{1}{\hat{\sigma}_{j(n)}}-\frac{1}
{\sigma_{j(n)}}\right|\right),\nonumber
\end{eqnarray}
where $x\vee y$ denotes the maximum of $x,y$ and
\begin{eqnarray}
\pi_{j(n)}&=&2\Phi\left(-\left|
Z_{j(n)}+\frac{\sqrt{n}(\mu_{j(n)}-\mu_{0,j(n)})}
{\sigma_{j(n)}}\right|\right).\label{c2.e2}
\end{eqnarray}
In particular, if $n\rightarrow\infty$, $\max_{1\leq j\leq p_n}
|\hat{\sigma}_{j(n)}-\sigma_{j(n)}|/(\sigma_{j(n)}\hat{\sigma}_{j(n)})
\rightarrow 0$ in probability, and
\begin{eqnarray}
\frac{\log(n\vee p_n)}{\sqrt{n}}\times
\max_{1\leq j\leq p_n}\frac{|b_{j(n)}-a_{j(n)}|}{\sigma_{j(n)}}
&\rightarrow& 0,\label{c2.e3}
\end{eqnarray}
then the left-hand-side of~(\ref{c2.e1}) $\rightarrow 0$ in probability.
\end{corollary}
\begin{remark}\label{c2.r1}
When $|b_{j(n)}-a_{j(n)}|/\sigma_{j(n)}$ is bounded, condition~(\ref{c2.e3})
becomes $\log^2 p_n/n=o(1)$. 
\end{remark}
\begin{remark}\label{r2}
Now, suppose the indices $J_n=\{1,\ldots,p_n\}$ are divided into two
groups, $J_{0n}$ and $J_{1n}$, where $H_0^{j(n)}$ holds for all 
$j\in J_{0n}$ and where
$\delta_{j(n)}=|\mu_{j(n)}-\mu_{0,j(n)}|/\sigma_{j(n)}>\tau$
for all $j\in J_{1n}$, where $\tau>0$. Then all of the
$\hat{\pi}_{j(n)}$s for $j\in J_{0n}$ will simultaneously converge to
uniform random variables with the same dependency structure inherent
in the data (as per the discussion before theorem~\ref{t1} above).  
Moreover, all of the
$\hat{\pi}_{j(n)}$ for $j\in J_{1n}$ will simultaneously converge to 0.
Thus the q-value approach to controlling FDR given in Storey, Taylor
and Siegmund (2004) should work under their weak dependence conditions
(7)--(9) (see also their theorem~5). A minor adjustment to this
argument will also work for contiguous alternative hypotheses where
the $\sqrt{n}\delta_{j(n)}$ quantities converge to bounded constants.
\end{remark}

\vspace{0.35cm}
{\it 2.3 Estimation of marginal sample medians.}
Now we consider inference for the median. Assume that each $F_{j(n)}$
has median $\xi_{j(n)}$ and is continuous in a neighborhood of
$\xi_{j(n)}$ with density $f_{j(n)}$. In this section, we do not
require the support of $F_{j(n)}$ to be compact. We do, however,
assume that there exists $\eta,\tau>0$ such that
\begin{eqnarray}
\min_{1\leq j\leq p_n}\inf_{x:|x-\xi_{j(n)}|\leq\eta}
f_{j(n)}(x)&\geq&\tau.\label{c4.e1}
\end{eqnarray}
Denote the sample median for $X_{1j(n)},\ldots,X_{nj(n)}$
as $\hat{\xi}_{j(n)}$.  More precisely, let $\hat{\xi}_{j(n)}=
\inf\{x:\hat{F}_{j(n)}(x)\geq 1/2\}$. The following corollary   
gives simultaneous consistency of $\hat{\xi}_{j(n)}$:
\begin{corollary}\label{c4}
Under condition~(\ref{c4.e1}) (for some $\eta,\tau>0$)
and the conditions of corollary~1, we have that
\begin{eqnarray}
\max_{1\leq j\leq p_n}|\hat{\xi}_{j(n)}-\xi_{j(n)}|&=&
O_P\left(\frac{\log(n\vee p_n)}{n}+\sqrt{\frac{\log p_n}{n}}\right).
\label{c4.e2}   
\end{eqnarray}
\end{corollary}

Now assume that we wish to test the marginal null hypothesis
$H_0^{j(n)}:\xi_{j(n)}=\xi_{0,j(n)}$ with the test statistics
$U_{j(n)}=2\sqrt{n}\hat{f}_{j(n)}(\hat{\xi}_{j(n)}-\xi_{0,j(n)})$,
where $\hat{f}_{j(n)}$ is a consistent estimator of $f_{j(n)}(\xi_{j(n)})$.
As duscussed in Kosorok~(1999), this is a good
choice of median test because it converges
rapidly to its limiting Gaussian distribution and appears to have
better moderate sample size performance compared to other median tests.
As with the marginal mean test, we need consistent estimates of the
p-values of these tests. We now study the consistency of the
p-value estimates $\hat{\pi}_{j(n)}'=2\Phi(-|U_{j(n)}|)$. We need   
some additional conditions.  Assume there exists $\eta,\tau>0$ and
$M<\infty$ such that~(\ref{c4.e1}) holds and, moreover, that
\begin{eqnarray}
\max_{1\leq j\leq p_n}\sup_{x:|x-\xi_{j(n)}|\leq\eta}f_{j(n)}&\leq&M
\label{c5.e1}
\end{eqnarray}  
and
\begin{eqnarray}
\max_{1\leq j\leq p_n}\sup_{\epsilon\leq\eta}\sup_{u:|u|\leq\epsilon}
\frac{|f_{j(n)}(\xi_{j(n)}+u)-f_{j(n)}(\xi_{j(n)})|}{\epsilon^{1/2}}
&\leq&M.\label{c5.e2}
\end{eqnarray}
We now have the following corollary:

\begin{corollary}\label{c5}
Under conditions~(\ref{c4.e1}), (\ref{c5.e1}) and~(\ref{c5.e2}),
for some $\eta,\tau>0$ and $M<\infty$, and provided both
$\max_{1\leq j\leq p_n}|\hat{f}_{j(n)}-f_{j(n)}(\xi_{j(n)})|=o_P(1)$
and $\log^3 p_n/n\rightarrow 0$ as $n\rightarrow\infty$, we have that
\begin{eqnarray}
\max_{1\leq j\leq p_n}|\hat{\pi}_{j(n)}'-\pi_{j(n)}'|&=&o_P(1),\label{c5.e3}
\end{eqnarray}
where
\begin{eqnarray}
\pi_{j(n)}'&=&2\Phi\left(-\left|Z_{j(n)}+2\sqrt{n}f_{j(n)}(\xi_{j(n)})
(\xi_{j(n)}-\xi_{0,j(n)})\right|\right),\label{c5.e3.n1}
\end{eqnarray}
and, for each $n\geq 1$, $Z_{1(n)},\ldots,Z_{p_n(n)}$ are
standard normals conditionally independent given ${\cal F}_n$ and for which
each $Z_{j(n)}$ has conditional distribution given ${\cal F}_n$
depending only on ${\cal F}_{j(n)}$, $1\leq j\leq p_n$.
\end{corollary}

Now, for corollary~\ref{c5} to be useful in conducting inference, we
need simultaneously consistent estimators $\hat{f}_{j(n)}$. One
possibility is
\begin{eqnarray}
\hat{f}_{j(n)}&=&\frac{\hat{F}_{j(n)}(\hat{\xi}_{j(n)}+\tilde{h}_{j(n)})
-\hat{F}_{j(n)}(\hat{\xi}_{j(n)}-\tilde{h}_{j(n)})}{2\tilde{h}_{j(n)}},
\label{kernel.e1}
\end{eqnarray}
where the window widths $\tilde{h}_{j(n)}$ are allowed to depend on the data
but must satisfy $\max_{1\leq j\leq p_n}\tilde{h}_{j(n)}=o_P(1)$ and
\begin{eqnarray}
\max_{1\leq j\leq p_n}\tilde{h}_{j(n)}^{-1}\left(\frac{\log n\vee p_n}{n}
+\sqrt{\frac{\log p_n}{n}}\right)&=&o_P(1).\label{kernel.e2}
\end{eqnarray}
If, in addition to the conditions of corollary~\ref{c5},
we assume conditions~(\ref{c4.e1}) and~(\ref{c5.e1}) apply to
the lower and upper quartiles of the distributions $F_{j(n)}$,
then $\tilde{h}_{j(n)}=2\hat{I}_{j(n)}n^{-1/5}$, where $\hat{I}_{j(n)}$
is the sample interquartile range based on $\hat{F}_{j(n)}$,
satisfies this requirement. This can be argued by first noting that
$\hat{I}_{j(n)}$ is asymptotically simultaneously bounded above
and below and that
\[n^{-1/5}\sqrt{\frac{\log p_n}{n}}=\sqrt{\frac{\log p_n}{n^{3/5}}} 
=\sqrt{\frac{\log p_n}{n^{1/3}}n^{-4/15}}\rightarrow 0.\]
There are many other possibilities that will also work.

\vspace{0.35cm}
{\bf 3. Marginal asymptotics for two-sample comparisons.}
The results of section~2 can be extended to two
sample results, where we have two i.i.d. samples of vectors of
length $p_n$, where $n=n_1+n_2$, and where $n_k$ is the size of
sample $k$, for $k=1,2$. Consistency results for estimating 
marginal distribution functions, marginal means and marginal
medians follows essentially without modification from theorem~\ref{t.c1} and
corollaries~\ref{c3} and~\ref{c4}.
Our interest will therefore focus on the more challenging issue
of testing whether the marginal means or medians are the same 
between the two samples.
We use superscript $(k)$ to denote membership in group $k$,
for $k=1,2$. In particular, $X_{i(n)}^{(k)}=(X_{i1(n)}^{(k)},
\ldots,X_{ip_n(n)}^{(k)})'$ is the $i$th observed vector in the 
$k$th group. In a similar manner, $F_{j(n)}^{(k)}$, ${\cal F}_{j(n)}^{(k)}$,
$a_{j(n)}^{(k)}\neq b_{j(n)}^{(k)}$, $\mu_{j(n)}^{(k)}$,
$\sigma_{j(n)}^{(k)}>0$, $\xi_{j(n)}^{(k)}$ and
$f_{j(n)}^{(k)}$, for $1\leq j\leq p_n$,
$k=1,2$, and all $n\geq 1$,
are the two-sample versions of the corresponding one-sample 
quantities introduced in section~2. Also let ${\cal F}_{j(n)}^{\ast}
=\sigma\left({\cal F}_{j(n)}^{(1)},{\cal F}_{j(n)}^{(2)}\right)$
and ${\cal F}_n^{\ast}=\sigma\left({\cal F}_{1(n)}^{\ast},
\ldots,{\cal F}_{p_n(n)}^{\ast}\right)$.

We first consider comparing the marginal means. Let 
$\bar{X}_{j(n)}^{(k)}$ be the sample mean of 
$X_{1j(n)}^{(k)},\ldots,$ $X_{nj(n)}^{(k)}$. Now assume that we wish
to test the marginal null hypothesis $H_0^{j(n)}:\mu_{j(n)}^{(1)}
=\mu_{j(n)}^{(2)}$ with the test statistic
\[T_{j(n)}^{\ast}=\sqrt{\frac{n_1 n_2}{n_1\left[\hat{\sigma}_{j(n)}^{(2)}
\right]^2+n_2\left[\hat{\sigma}_{j(n)}^{(1)}\right]^2}}
\left(\bar{X}_{j(n)}^{(1)}-\bar{X}_{j(n)}^{(2)}\right),\]
where $\hat{\sigma}_{j(n)}^{(k)}$ is a location-invariant
and consistent estimator of $\sigma_{j(n)}^{(k)}$, $k=1,2$. The
following corollary provides conditions under which p-values
estimated by $\hat{\pi}_{j(n)}^{\ast}=2\Phi\left(-|T_{j(n)}^{\ast}|\right)$
are uniformly consistent over all $1\leq j\leq p_n$:
\begin{corollary}\label{c6}
Let the constants $c_1,c_2$ be as in theorem~\ref{t1}. Then for
all $n_1,n_2,p_n\geq 2$, there exist standard normal random variables
$Z_{1(n)}^{\ast},\ldots,Z_{p_n(n)}^{\ast}$ which are conditionally
independent given ${\cal F}_n^{\ast}$ and for which each
$Z_{j(n)}^{\ast}$ has conditional distribution given ${\cal F}_n^{\ast}$
depending only on ${\cal F}_{j(n)}^{\ast}$, $1\leq j\leq p_n$, such that
\begin{eqnarray}
\;\;\;\;\max_{1\leq j\leq p_n}
\left|\hat{\pi}_{j(n)}^{\ast}-\pi_{j(n)}^{\ast}\right|&\leq&
\sum_{k=1,2}\left[
\frac{c_1\log n_k+c_2\log p_n}{\sqrt{n_k}}\left(\max_{1\leq j\leq p_n}
\frac{|b_{j(n)}^{(k)}-a_{j(n)}^{(k)}|}
{\sigma_{j(n)}^{(k)}}\right)\right.\label{c6.e1}\\  
&&\left.+\frac{1}{2}\left(\max_{1\leq j\leq n}
\left(\hat{\sigma}_{j(n)}^{(k)}\vee\sigma_{j(n)}^{(k)}\right)
\left|\frac{1}{\hat{\sigma}_{j(n)}^{(k)}}-\frac{1}
{\sigma_{j(n)}^{(k)}}\right|\right)\right],\nonumber
\end{eqnarray}
where 
\begin{eqnarray}
\pi_{j(n)}^{\ast}&=&2\Phi\left(-\left|
Z_{j(n)}^{\ast}+\sqrt{\frac{n_1 n_2}{n_1\left[\sigma_{j(n)}^{(2)}\right]^2
+n_2\left[\sigma_{j(n)}^{(1)}\right]^2}}
\left(\mu_{j(n)}^{(1)}-\mu_{j(n)}^{(2)}\right)\right|\right).
\label{c6.e2}
\end{eqnarray}
In particular, if $n_k\rightarrow\infty$, $\max_{1\leq j\leq p_n}
\left|\hat{\sigma}_{j(n)}^{(k)}-\sigma_{j(n)}^{(k)}\right|/
\left(\sigma_{j(n)}^{(k)}\hat{\sigma}_{j(n)}^{(k)}\right)
\rightarrow 0$ in probability, and
\begin{eqnarray}
\frac{\log(n_k\vee p_n)}{\sqrt{n_k}}\times
\max_{1\leq j\leq p_n}\frac{\left|b_{j(n)}^{(k)}-a_{j(n)}^{(k)}\right|}
{\sigma_{j(n)}^{(k)}}&\rightarrow& 0,\label{c6.e3}
\end{eqnarray}
for $k=1,2$, then the left-hand-side of~(\ref{c6.e1}) 
$\rightarrow 0$ in probability.
\end{corollary}

We now consider comparing marginal medians. Assume that we wish to
test the marginal null hypothesis $H_0^{j(n)}:\xi_{j(n)}^{(1)}
=\xi_{j(n)}^{(2)}$ with the test statistic
\[U_{j(n)}^{\#}=2\sqrt{\frac{n_1 n_2}{n_1/\left[\hat{f}_{j(n)}^{(2)}
\right]^2+n_2/\left[\hat{f}_{j(n)}^{(1)}\right]^2}}\left(
\hat{\xi}_{j(n)}^{(1)}-\hat{\xi}_{j(n)}^{(2)}\right),\]
where $\hat{f}_{j(n)}^{(k)}$ is consistent for $f_{j(n)}^{(k)}
(\xi_{j(n)}^{(k)})$, $k=1,2$. The following corollary provides conditions
under which p-values estimated by $\hat{\pi}^{\#}_{j(n)}
=2\Phi\left(-|U_{j(n)}^{\#}|\right)$ are uniformly consistent
over all $1\leq j\leq p_n$:
\begin{corollary}\label{c7}
Assume that the one-sample
conditions given in expressions~(\ref{c4.e1}), (\ref{c5.e1})
and~(\ref{c5.e2}), for all of the marginal distribution functions and 
densities in both samples, are satisfied for constants $\eta,\tau>0$
and $0<M<\infty$. Assume also that
$\max_{\,k=1,2;\,1\leq j\leq p_n}\left|\hat{f}_{j(n)}^{(k)}-
f_{j(n)}^{(k)}(\xi_{j(n)}^{(k)})\right|=o_P(1)$ and
$\log^3 p_n/(n_1\wedge n_2)\rightarrow 0$ as
$n\rightarrow\infty$. Then
\begin{eqnarray}
\max_{1\leq j\leq p_n}\left|\hat{\pi}_{j(n)}^{\#}-\pi_{j(n)}^{\#}
\right|=o_P(1),\label{c7.e1}
\end{eqnarray}
where
\begin{eqnarray}
\label{c7.e2}&&\\
\pi_{j(n)}^{\#}=2\Phi\left(-\left|Z_{j(n)}^{\ast}+
2\sqrt{\frac{n_1 n_2}{n_1/\left[f_{j(n)}^{(2)}(\xi_{j(n)}^{(2)})\right]^2
+n_2/\left[f_{j(n)}^{(1)}(\xi_{j(n)}^{(1)})\right]^2}}
\left(\xi_{j(n)}^{(1)}-\xi_{j(n)}^{(2)}\right)\right|\right),&&\nonumber
\end{eqnarray}
and, for each $n\geq 1$, $Z_{1(n)}^{\ast},\ldots,Z_{p_n(n)}^{\ast}$ are
standard normals conditionally independent given ${\cal F}_{n}^{\ast}$
and for which each $Z_{j(n)}^{\ast}$ has conditional distribution
given ${\cal F}_n^{\ast}$ depending only on ${\cal F}_{j(n)}^{\ast}$,
$1\leq j\leq p_n$.
\end{corollary}

\vspace{0.35cm}
{\bf 4. Distribution free statistics.} When the distribution of the
test statistic under the null hypothesis does not depend on the
distribution function, results stronger than those presented in
sections~2 and~3 are possible for marginal p-value consistency.
Consider first the one-sample setting, and assume that the
distributions $F_{j(n)}$ are all continuous and symmetric around
their respective medians. Suppose we are interested in marginal
testing of $H_0^{j(n)}:\xi_{j(n)}=0$ 
using the signed rank test $\tilde{T}_{j(n)}$ studied in section~3 of
Kosorok and Ma (2005). Define 
\[V_{j(n)}=\frac{\tilde{T}_{j(n)}-(n^2+n)/4}{\sqrt{(3 n^3+2 n^2+n)/24}}.\]
Note that the distribution of $V_{j(n)}$ does not depend on
$F_{j(n)}$ under $H_0^{j(n)}$. 
Let $\Phi_n$ be the exact distribution of $V_{j(n)}$
under $H_0^{j(n)}$. It is easy to verify that $\Phi_n$ converges 
uniformly to $\Phi$. Hence   
\[\max_{1\leq j\leq p_n}|2\Phi_n(-|V_{j(n)}|)-2\Phi(-|V_{j(n)}|)|
\rightarrow 0,\]
regardless of how fast $p_n$ grows. Thus the normal approximation
is simultaneously consistent for the true p-values 
when $n\rightarrow\infty$, without any constraints on $p_n$.

The key feature that makes this work is that the p-values depend 
only on the correctness of the probability calculation under the null
hypothesis. P-value computations do not require knowledge of the
distribution under alternatives. The only possibly unnatural assumption
required for the above signed-rank test is symmetry about the median.
An alternative statistic is the sign test.  Under the null hypothesis
that the median is zero, the sign test is Bernoulli with probability
$1/2$. As with the signed-rank test, the standardized
sign test under the null converges to a normal limit.  A disadvantage
of the sign test is that the range of possible values is limited,
resulting in a granular distribution which converges somewhat
slowly to the normal limit. 

Similar reasoning applies to distribution-free two-sample test
statistics.  Interestingly, there appears to be a larger variety
of useful tests to choose from which do not require specification of the
distribution function than there are in the one-sample
setting. Suppose we are interested
in marginal testing of $H_0^{j(n)}:F_{j(n)}^{(1)}=F_{j(n)}^{(2)}$,
and we assume that the $F_{j(n)}^{(k)}$ are continuous
for all $1\leq j\leq p_n$ and $k=1,2$. Let
$\hat{F}_{j(n)}^{(k)}(t)=n_k^{-1}\sum_{i=1}^{n_k}\ind\left\{
X_{ij(n)}^{(k)}\leq t\right\}$, for $k=1,2$; 
$\hat{F}_{j(n)}^{(0)}=n^{-1}\left[n_1\hat{F}_{j(n)}^{(1)}
+n_2\hat{F}_{j(n)}^{(2)}\right]$; and
$\hat{G}_{j(n)}=\sqrt{n_1n_2/n}\left(\hat{F}_{j(n)}^{(1)}
-\hat{F}_{j(n)}^{(2)}\right)$. We now consider several
statistics which are invariant under monotone transformations of
the data:
\begin{enumerate}
\item The two-sample Wilcoxon rank sum test 
$\tilde{T}_{j(n)}^{\ast 1}=\sqrt{12}
\int_{\re}\hat{G}_{j(n)}(s)d\hat{F}_{j(n)}^{(0)}(s)$;
\item The two-sample Kolmogorov-Smirnov test 
$\tilde{T}_{j(n)}^{\ast 2}=\sup_{t\in\re}
\left|\hat{G}_{j(n)}\right|$;
\item The two-sample Cram\'{e}r-von Mises test
$\tilde{T}_{j(n)}^{\ast 3}=\int_{\re}\hat{G}_{j(n)}^2(s)
d\hat{F}_{j(n)}^{(0)}(s)$.
\end{enumerate}

Fix $j\in\{1,\ldots,p_n\}$ and
assume $H_0^{j(n)}$ holds. All three of these statistics
are now invariant under the monotone transformation
$t\mapsto F_{j(n)}^{(0)}(t)$, where $F_{j(n)}^{(0)}\equiv
F_{j(n)}^{(1)}=F_{j(n)}^{(2)}$. Thus, without loss of
generality, we can assume the data are i.i.d. uniform~$[0,1]$.
For $m=1,2,3$,
let $K_n^{\ast m}$ be the corresponding cumulative distribution 
function for the statistic $\tilde{T}_{j(n)}^{\ast m}$ under this
uniformity assumption (note that it does not depend on $j$
because of the invariance), and let $K_0^{\ast m}$ be the limiting
cumulative distribution function. Suppose that we compute 
approximate p-values for the three statistics as follows: 
$\hat{\pi}_{j(n)}^{\ast 1}=2\Phi\left(-\left|\tilde{T}_{j(n)}^{\ast 1}
\right|\right)$, and $\hat{\pi}_{j(n)}^{\ast m}=
1-K_0^{\ast m}\left(\tilde{T}_{j(n)}^{\ast m}\right)$, for $m=2,3$.
Because it can be shown that $K_0^{\ast m}$ is continuous for all
$m=1,2,3$, the convergence of $K_n^{\ast m}$ to $K_0^{\ast m}$ is
uniform. Thus, even after we drop the $H_0^{j(n)}$ assumption,
the approximate p-values based on $K_n^{\ast m}$
are simultaneously consistent for the true p-values, as
$n\rightarrow\infty$, without constrainting $p_n$.

The following lemma yields the form of $K_0^{\ast m}$, for $m=1,2,3$.
The results are essentially
classical, but they are included here for completeness:
\begin{lemma}\label{l2}
For $m=1,2,3$, $K_n^{\ast m}$ converges uniformly to $K_0^{\ast m}$, 
as $n_1\wedge n_2\rightarrow\infty$, where
\begin{itemize}
\item $K_0^{\ast 1}=\Phi$;
\item For $t>0$, $K_0^{\ast 2}(t)=1-2\sum_{l=1}^{\infty}(-1)^le^{-2l^2t^2}$
is the distribution of the supremum in absolute value of
a standard Brownian bridge;
\item $K_0^{\ast 3}$ is the distribution of $\pi^{-2}\sum_{l=1}^{\infty}
l^{-2}\tilde{Z}_l^2$, where $\tilde{Z}_1,\tilde{Z}_2,\ldots$ are
i.i.d. standard normals.
\end{itemize}
\end{lemma}

\vspace{0.35cm}
{\bf 5. Impact of microarray normalization.}
In this section, we consider the affect of normalization on the 
theory presented in sections~2--4. 
For the simple normalization
model~(\ref{model.1}), this will require the $\hat{h}_i$s to be uniformly
consistent at the rate $O_P(\sqrt{n}\log n)$. This requirement seems
reasonable for certain estimation methods, including the method described
in FPH.  In this method, data across all genes within each array are
utilized for estimating the $h_i$s. Since the
number of genes $p_n$ usually increases nearly exponentially relative
to the number of microarrays, the
number of observations available for estimating the $h_i$s is
many orders of magnitude higher than $n$, even after taking into
account dependencies within arrays and the fact that the number
of arrays is increasing in~$n$. For this particular facet of our
problem, the large number of genes actually works in our favor.
A variant of this argument can also be found
in Kosorok and Ma (2005).

Consider first the one-sample setting of section~2. Let $\tilde{X}_{i(n)}
=(\tilde{X}_{i1(n)},\ldots,\tilde{X}_{ip_n(n)})'$ be an approximation
of the ``true data'' $X_{i(n)}$, $1\leq i\leq n$, 
and define 
\[\hat{\epsilon}_n=
\max_{1\leq j\leq p_n;\,1\leq i\leq n}|\tilde{X}_{ij(n)}-X_{ij(n)}|.\]
With proper, partially consistent normalization, the true gene
effects $\{X_{ij(n)},1\leq j\leq p_n,1\leq i\leq n\}$ should be
uniformly consistently estimated by the residuals from the
normalization $\{\tilde{X}_{ij(n)},1\leq j\leq p_n,1\leq i\leq n\}$.
In other words $\hat{\epsilon}_n=o_P(1)$. The essence of our arguments
involves an assessment of how well $\tilde{F}_{j(n)}(t)
\equiv n^{-1}\sum_{i=1}^n\ind\{\tilde{X}_{ij(n)}\leq t\}$ approximates
$\hat{F}_{j(n)}(t)$ uniformly in $t$. We need the following 
strengthening of condition~(\ref{c5.e1}):
\begin{eqnarray}
\limsup_{n\rightarrow\infty}
\max_{1\leq j\leq p_n}\sup_{t\in\re}f_{j(n)}(t)&\leq&\tilde{M},\label{t3.e0}
\end{eqnarray}
for some $\tilde{M}<\infty$. We now have the following theorem, the
proof of which involves a precise bound on the modulus of continuity
of Brownian motion (see lemma~\ref{l1} in section~8 below):
\begin{theorem}\label{t3}
Assume condition~(\ref{t3.e0}) holds for some $\tilde{M}<\infty$. Then
the following are true:
\begin{enumerate} 
\item[(i)] If $\log p_n/n=o(1)$ and $\hat{\epsilon}_n=o_P(1)$, then
\[\max_{1\leq j\leq p_n}
\left\|\tilde{F}_{j(n)}-\hat{F}_{j(n)}\right\|_{\infty}=o_P(1);\]
\item[(ii)] If, in addition, $\log^2 p_n/n=o(1)$ and 
$\sqrt{n}(\log n)\hat{\epsilon}_n=O_P(1)$, then also
\[\max_{1\leq j\leq p_n}
\left\|\tilde{F}_{j(n)}-\hat{F}_{j(n)}\right\|_{\infty}=o_P(n^{-1/2}).\]
\end{enumerate}
\end{theorem}

\begin{remark}\label{t3.r1}
Note that the one-sample signed rank test $\tilde{T}_{j(n)}$ can be written
as a normalization of $\sqrt{n}\int_{\re}\left[\hat{F}_{j(n)}(u)
-\hat{F}_{j(n)}(-u)\right]d\hat{F}_{j(n)}(u)$, and the one-sample
sign test can be written as a normalization of
$\sqrt{n}\int_{\re}\mbox{sign}(u)d\hat{F}_{j(n)}(u)$. Thus part~(ii)
of theorem~\ref{t3}
allows us to replace $\hat{F}_{j(n)}$ with $\tilde{F}_{j(n)}$ in
both of these statistics without destroying the simultaneous
consistency over $1\leq j\leq p_n$ established in section~4
of the normal approximation for the true p-values based on
the true data. 
\end{remark}

Theorem~\ref{t3} can also be used to verify that the 
asymptotic results for the one-sample mean and median tests of sections~2
and~3 can be similarly extended for the approximate data 
$\tilde{X}_{1(n)},\ldots,\tilde{X}_{n(n)}$. For $j=1,\ldots,p_n$, let
$\check{X}_{j(n)}$ be the sample
mean of $\tilde{X}_{1j(n)},\ldots,\tilde{X}_{nj(n)}$,
and define the approximate sample median
$\tilde{\xi}_{j(n)}=\inf\{r:\tilde{F}_{j(n)}(r)\geq 1/2\}$. The
following corollary yields consistency of these estimators:
\begin{corollary}\label{c8}
Assume the conditions of theorem~\ref{t3}, part~(i), hold. Then
\begin{enumerate}
\item[(i)] $\max_{1\leq j\leq p_n}\left\|\tilde{F}_{j(n)}-F_{j(n)}
\right\|_{\infty}=o_P(1)$;
\item[(ii)] Provided $\limsup_{n\rightarrow\infty}\max_{1\leq j\leq p_n}
|b_{j(n)}-a_{j(n)}|<\infty$, $\max_{1\leq j\leq p_n}|\check{X}_{j(n)}
-\mu_{j(n)}|=o_P(1)$;
\item[(iii)] $\max_{1\leq j\leq p_n}|\tilde{\xi}_{j(n)}-\xi_{j(n)}|=o_P(1)$.
\end{enumerate}
\end{corollary}

The following corollary strengthens result~(ii) of corollary~\ref{c8} above
and yields consistency of the p-values of one-sample tests based on the 
approximate data:
\begin{corollary}\label{c9}
Assume the conditions of theorem~\ref{t3}, part~(ii), hold. Then the
following results are true under the given conditions:
\begin{enumerate}
\item[(i)] Provided $\limsup_{n\rightarrow\infty}\max_{1\leq j\leq p_n}
n^{1/4}|b_{j(n)}-a_{j(n)}|<\infty$, $\max_{1\leq j\leq p_n}|\check{X}_{j(n)}
-\mu_{j(n)}|=o_P(1)$.
\item[(ii)] Suppose that the conditions of corollary~\ref{c2} hold,
except that $\check{X}_{j(n)}$ is used instead of $\bar{X}_{j(n)}$
and that all other estimated quantities are based on $\tilde{F}_{j(n)}$
rather than on $\hat{F}_{j(n)}$, for $j=1,\ldots,p_n$. Then, provided
\[\limsup_{n\rightarrow\infty}\max_{1\leq j\leq p_n}\frac{\left\{
|b_{j(n)}-a_{j(n)}|\vee 1\right\}}{\sigma_{j(n)}}<\infty\] 
and $\max_{1\leq j\leq p_n}\left|
\hat{\sigma}_{j(n)}-\sigma_{j(n)}\right|/\left(
\sigma_{j(n)}\hat{\sigma}_{j(n)}\right)\rightarrow 0$, 
$\max_{1\leq j\leq p_n}\left|\hat{\pi}_{j(n)}-\pi_{j(n)}
\right|=o_P(1)$, for the filtrations ${\cal F}_n$ and ${\cal F}_{j(n)}$, 
$1\leq j\leq p_n$, based on the true data.
\item[(iii)] Suppose that the conditions of corollary~\ref{c5} hold, except
that $\tilde{\xi}_{j(n)}$ is used instead of $\hat{\xi}_{j(n)}$ 
and that all other estimated quantities are based on $\tilde{F}_{j(n)}$
rather than on $\hat{F}_{j(n)}$, for $j=1,\ldots,p_n$. Then the conclusions
of corollary~\ref{c5} still hold for the filtrations ${\cal F}_n$
and ${\cal F}_{j(n)}$, $1\leq j\leq p_n$, based on the true data.
\end{enumerate}
\end{corollary}

\begin{remark}\label{c9.r1}
Parts~(ii) and~(iii) of corollary~\ref{c9} tell us that we can construct
valid mean and median based hypothesis tests from suitably normalized data, 
and that any dependencies beyond the original dependency structure
induced by the approximation vanish asymptotically. Thus the arguments given
in remark~\ref{r2} regarding the validity of the q-value approach
for controlling FDR still hold after normalization.
\end{remark}

The extension of these results to the two-sample setting is straightforward.
As done in section~4, we will use superscript~$(k)$ to denote membership
in group $k$, for $k=1,2$. Let $\tilde{F}_{j(n)}^{(k)}$ be the empirical 
distribution of the approximate data sample $\tilde{X}_{1j(n)}^{(k)},\ldots,
\tilde{X}_{nj(n)}^{(k)}$; $\tilde{\xi}_{j(n)}^{(k)}=\inf\{r:
\tilde{F}_{j(n)}(r)\geq 1/2\}$; $\hat{\epsilon}_n^{(k)}$ be the maximum
error between the approximate and true data for group $k$; 
and redefine $\hat{\epsilon}_n=\hat{\epsilon}_n^{(1)}\vee
\hat{\epsilon}_n^{(2)}$. Also let $\check{T}_{j(n)}^{\ast m}$ be the
version of $\tilde{T}_{j(n)}^{\ast m}$ with $\tilde{F}_{j(n)}^{(k)}$
replacing $\hat{F}_{j(n)}^{(k)}$, for $k=1,2$ and $m=1,2,3$.
The following corollary gives the main two-sample approximation results:
\begin{corollary}\label{c10}
Assume $n_1\wedge n_2\rightarrow\infty$;
$\limsup_{n\rightarrow\infty}\max_{k=1,2}\max_{1\leq j\leq p_n}
\|f_{j(n)}^{(k)}\|_{\infty}\leq\tilde{M}$, for some $\tilde{M}<\infty$;
$\log^2 p_n/(n_1\wedge n_2)=o(1)$; and  
$\sqrt{n}(\log n)\hat{\epsilon}_n=O_P(1)$. Then the following are true
under the given conditions:
\begin{enumerate}
\item[(i)] Suppose that the conditions of corollary~\ref{c6} hold,
except the sample means are based on the approximate data
and all other estimated quantities are based on $\tilde{F}_{j(n)}^{(k)}$
rather than on $\hat{F}_{j(n)}^{(k)}$, for $j=1,\ldots,p_n$ and $k=1,2$. 
Then, provided
\[\limsup_{n\rightarrow\infty}\max_{k=1,2}\max_{1\leq j\leq p_n}
\frac{\left\{|b_{j(n)}^{(k)}-a_{j(n)}^{(k)}|\vee 1\right\}}
{\sigma_{j(n)}^{(k)}}<\infty\] 
and $\max_{k=1,2}\max_{1\leq j\leq p_n}\left|
\hat{\sigma}_{j(n)}^{(k)}-\sigma_{j(n)}^{(k)}\right|/\left(
\sigma_{j(n)}^{(k)}\hat{\sigma}_{j(n)}^{(k)}\right)\rightarrow 0$, 
\[\max_{1\leq j\leq p_n}\left|\hat{\pi}_{j(n)}^{\ast}-\pi_{j(n)}^{\ast}
\right|=o_P(1),\] 
for the filtrations ${\cal F}_n$ and ${\cal F}_{j(n)}$, 
$1\leq j\leq p_n$, based on the true data.
\item[(ii)] Suppose that the conditions of corollary~\ref{c7} hold, except
that $\tilde{\xi}_{j(n)}^{(k)}$ is used instead of $\hat{\xi}_{j(n)}^{(k)}$ 
and that all other estimated quantities are based on $\tilde{F}_{j(n)}^{(k)}$
rather than on $\hat{F}_{j(n)}^{(k)}$, for $j=1,\ldots,p_n$ and $k=1,2$. 
Then the conclusions
of corollary~\ref{c7} still hold for the filtrations ${\cal F}_n$
and ${\cal F}_{j(n)}$, $1\leq j\leq p_n$, based on the true data.
\item[(iii)] $\max_{1\leq j\leq p_n}\left|
\check{T}_{j(n)}^{\ast m}-\tilde{T}_{j(n)}^{\ast m}\right|=o_P(1)$, 
for $m=1,2,3$. Thus the approximate p-values
based on the approximate data for the three distribution-free two-sample
tests given in section~4 are uniformly consistent for the true p-values
based on the true data.
\end{enumerate}
\end{corollary}

\vspace{0.35cm}
{\bf 6. Numerical studies.} 

{\it 6.1 One-sample simulation study.} We used a small simulation study to
assess the finite sample performance of the following one-sample
methodologies: (1) the mean based comparison of section 2.2, (2)
the median based comparison of section 2.3
and (3) the signed rank test of section 4. We set the number of genes
to $p=2000$ and the number of arrays to $n=20,50$. 
Let $Z_{i1},Z_{i2},\ldots$, $i=1,\ldots,n$, be a sequence of i.i.d.
standard normal random variables. We generated simulated data using
the following three models:

Model 1: $X_{ij}=H(Z_{ij})$ for $i=1, \ldots, p$;

Model 2: $X_{ij}=H\left (\sum_{l=(j-1)\times m+1}^{(j-1)\times m +k} 
Z_{l}/\sqrt{k}\right)$ with $k=10, m=7$;

Model 3: Same as Model 2, but with $k=10, m=3$.

In the above, $H=2 \Phi -1$, where $\Phi$ is the
cumulative distribution for the standard normal. This yields a
marginal $unif[-1,1]$ distribution for all three models. The
genes in model~1 are i.i.d., while in model~2 there is
strong dependence and in model~3 weak dependence between genes.
We assume the first 40 genes
have non-zero means, denoted as $\beta$ and generated from $unif[-2, 2]$.
For each approach, marginal p-values are computed based on the asymptotic
results for one-sample tests given in sections~2 and~4.  
For the median approach, density estimation is based on the interquartile
range band-width kernel described in the last paragraph of section~2.3.
We employ standard FDR techniques with expected FDR
$E(FDR)=0.2$. The marginal p-values are ranked, resulting in the
ordered p-values $\pi_{(1)}\leq \pi_{(2)}\leq \ldots\leq \pi_{(p)}$. 
Let $\tilde{g}$ be the largest $g$ such that $\pi_{(g)}\leq g/p \times q$,
where $q$ is the target FDR (for the simulations, $q=0.2$). 
Genes corresponding to $\pi_{(1)}, \ldots \pi_{(\tilde{g})}$ are identified 
as significantly differentially expressed.

Simulation results based on 100 replicates per scenario are shown in Table~1.
We can see that as the sample size increases, the performances of
all three approaches generally improve. When the sample size is small,
the mean based approach can effectively identify differentially
expressed genes, but with high false positive rates. Empirical FDRs 
for the rank approach are quite low. The
rank based approach misses quite a few true positives.
When the sample size is large, the median approach and the rank approach  
perform much better than the mean based approach, with less false
positives while still being able to identify true positives. The presence of
correlation appears to have very little impact on the performance.

\vspace{0.3cm}
{\it 6.2 Two-sample simulation study}.
Since the affect of dependence between genes in the simulation study 
of section~6.1 was minimal, we decided to restrict our focus on the
i.i.d. gene setting for the two-sample simulations. 
We set the number of genes to $p=2000$ and numbers of arrays
(sample sizes) to $n_1=n_2=10,30,60$. The model we explore
is Model 4: $X_{ij}^{(k)}\sim unif[-1, 1]$, $i=1,\ldots,k_k$,
$j=1,\ldots,p$, and $k=1,2$. For this data,
we apply the mean approach, the median approach, the Wilcoxon test and 
the Kolmogorov-Smirnov test to the two-sample comparison of
$X_{ij}^{(1)}+\beta$ versus $X_{ij}^{(2)}$, where $\beta$ is generated 
as in section~6.1 for the first 40 genes of each array.
Summary statistics for $E(FDR)=0.2$ and 100 replicates are shown in 
Table~2. Similar conclusions as in section~6.1
on the effects of sample size and gene distribution can be made. 
We especially notice that when
the sample size is small, the mean based approach appears to be the only one
that can identify a significant number of true positives.
The false positive rates are smaller than the
target for the median, Wilcoxon and Kolmogorov-Smirnov (KS) approaches. 
The mismatch between 
the empirical FDR with the target FDR can be serious for the mean approach,
especially when the sample size is small.

Based on other numerical studies (not presented), it appears that part
of the convergence difficulties with the nonparametric
approaches (in both the one and two sample settings) 
are due to the small number of distinct possible values
these statistics can have. 
It is unclear how to solve this problem for the nonparametric
one-sample tests, but
it appears that the two-sample tests can be improved by replacing
$\hat{G}_{j(n)}$ with $\check{G}_{j(n)}=\sqrt{n_1 n_2/n}
\left[\hat{F}_{j(n)}^{(1)}-\check{F}_{j(n)}^{(2)}\right]$, 
where $\check{F}_{j(n)}^{(2)}=(n_2/(n_2+1))\hat{F}_{j(n)}^{(2)}$.
This increases the number of possible values of the statistic, and 
preliminary simulation studies (also not presented) indicate that the rate 
of convergence for smaller sample sizes is improved. Thus we recommend
that this modification be considered whenever $n_1=n_2$.
Note that the modification does not affect the asymptotics since
\[\sqrt{\frac{n_1 n_2}{n}}\left|1-\frac{n_2}{n_2+1}\right|
\leq\frac{1}{\sqrt{n_2}}.\]

\vspace{0.3cm}
{\it 6.3 Estrogen data}. These datasets were first presented by West et al.
(2001) and Spang et al. (2001). Their common expression
matrix monitors 7129 genes in 49 breast tumor samples. The data
were obtained by applying the Affymetrix gene chip technology. 
The response describes the lymph nodal (LN) status, which is an
indicator for the metastatic spread of the tumor, an important
risk factor for disease outcome. 25 samples are positive (LN+) and  
24 samples are negative (LN-).
The goal is to identify genes differentially expressed between positive
and negative samples from the 3332 genes passing the first step 
of processing described in Dudoit, Fridlyand and Speed (2002). A base~2 
logarithmic transformation of the gene expressions is first applied.

\begin{figure}[htbp]
  \begin{center}
     \scalebox{0.50}[0.9]{\rotatebox{270}{\includegraphics{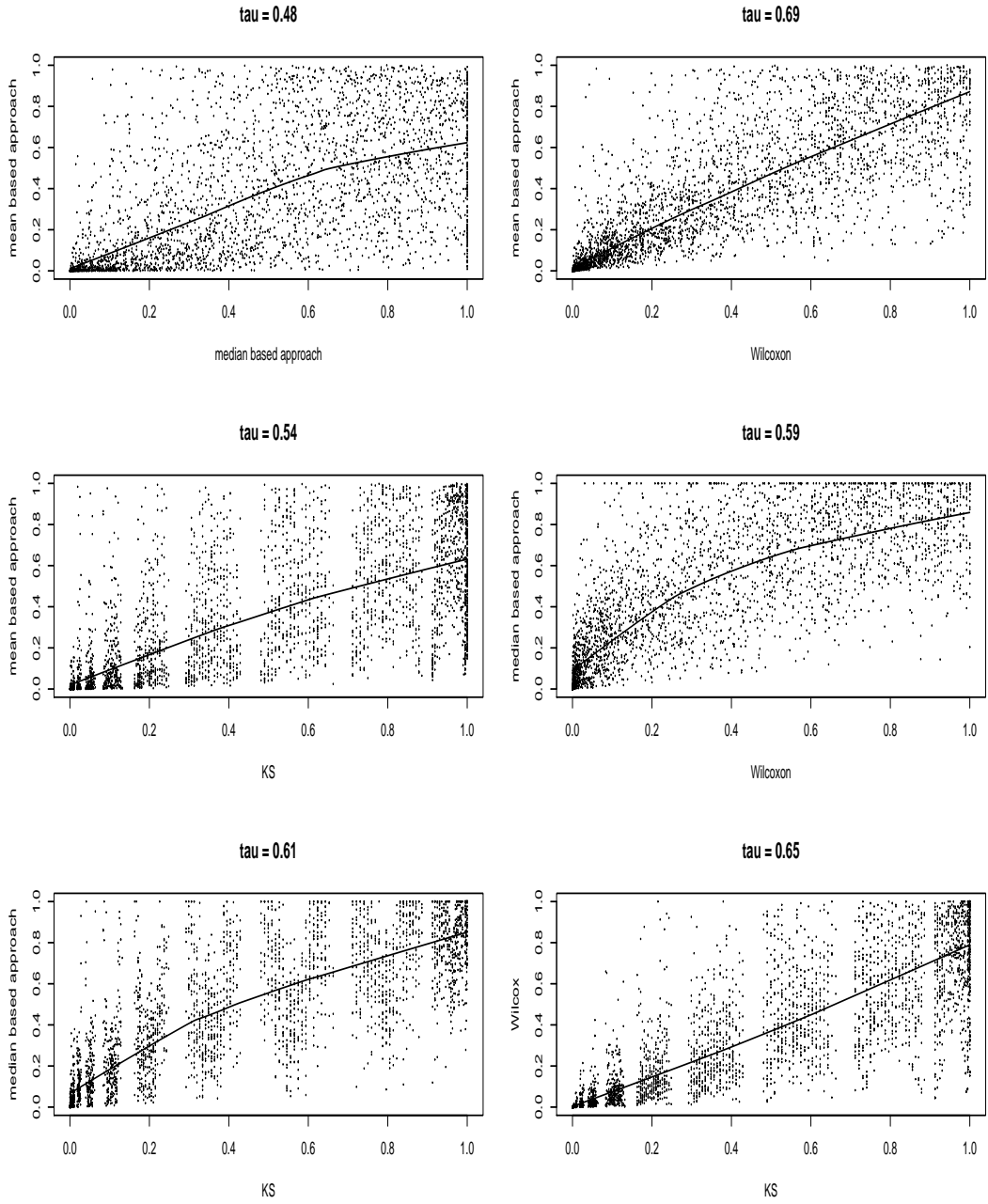}}}
        \caption{Estrogen data. Scatter plots of p-values comparing the four 
    approaches (mean, median, Wilcoxon and KS). A lowess smoother is used
    to estimate the trend, and the associated rank correlation coefficient 
    (tau) is given above each panel.}
    \label{fig:estrogen}
  \end{center}
\end{figure}

We set the target FDR to 0.1 and apply the standard FDR method with the four 
two-sample comparison approaches: 445 (mean), 261 (median), 423 (Wilcox)
and 211 (KS) genes are identified, respectively. The mean based approach 
and the Wilcoxon test identify significantly more
genes than the median approach and the KS test. This pattern was also
demonstrated in Table~2 (for sample size $n_1=n_2=30$).
It is unclear what causes these differences.
However, the overlaps of genes 
identified by the different approaches are substantial. 
For example, there are 196 common genes between the mean approach and the 
median approach. 
In Figure~\ref{fig:estrogen}, we show scatter plots of p-values from 
the different approaches. The rank correlation coefficients show
substantial similarities among different approaches. Note the banded
pattern in the plots involving the KS statistic. This is a consequence
of the low number of distinct possible values this statistic can have
as was discussed in section~6.2 above.

\vspace{0.35cm}
{\bf 7. Discussion.} 
The main results of this paper are that marginal 
(gene specific) estimates and asymptotic-based p-values are uniformly 
consistent in microarray experiments with $n$ replications---regardless
of the dependencies between genes---provided the 
number of genes $p_n$ satisfies $\log p_n=o(n)$, 
$\log p_n=o(n^{1/2})$ or $\log p_n=o(n^{1/3})$, depending on the desired task.
In other words, the number of genes is allowed to increase almost 
exponentially fast relative to the number of arrays. This seems to be a 
realistic asymptotic regime for microarray studies. These results also hold 
true for two-sample comparisons. Moreover, the results continue to hold even 
after normalization, provided the normalization process is sufficiently 
accurate. 

We note that the simulation and data analyses seem to 
support the theoretical results of the paper, although some test procedures 
appear to work better than others. We also acknowledge that a number of 
important issues, such 
as the affect of marginal distribution on the asymptotics and the affect of
normalization, were not evaluated in the limited simulation studies 
presented in section~6. A refined and more thorough simulation study that 
addresses these points is beyond the scope of the current paper but is worth 
pursuing in the future.

A theoretical limitation of the present study is that the asymptotics 
developed are not yet accurate enough to provide precise guidelines on 
sample size for specific microarray experiments. The development of such 
guidelines is worthwhile to pursue as a future topic, but it most likely
would require at least some assumptions on the dependencies between genes. 
Such assumptions are out of place in the present paper since a
strength of the paper is the absence of assumptions on gene
interdependence. It is because of this generality that we believe the 
results of this paper should be a useful point of departure for future, 
more refined asymptotic analyses of microarray experiments.

\newpage
{\bf 8. Proofs.}

{\it Proof of theorem~\ref{t.c1}.}
Define $V_{j(n)}\equiv\sqrt{n}\|\hat{F}_{j(n)}-F_{j(n)}\|_{\infty}$, and
note that by theorem~\ref{t.c1.t1} below combined with lemma~2.2.1 of
van der Vaart and Wellner (1996) (abbreviated VW hereafter),
$\|V_{j(n)}\|_{\psi_2}\leq \sqrt{3/2}$ for all $1\leq j\leq p_n$. 
Now, by lemma~2.2.2 of VW combined with the fact that
$\limsup_{x,y\rightarrow\infty}\psi_2(x)\psi_2(y)/\psi_2(xy)=0$,
we have that there exists a universal constant $c_{\ast}<\infty$
such that $\left\|\max_{1\leq j\leq p_n}V_{j(n)}\right\|_{\psi_2}
\leq c_{\ast}\sqrt{\log(1+p_n)}\sqrt{3/2}$ for all $n\geq 1$. The desired 
result now follows for the constant $c_0=\sqrt{6}c_{\ast}$,
since $\log(k+1)\leq 2\log k$ for any $k\geq 2$.$\Box$

\begin{theorem}\label{t.c1.t1}
Let $Y_1,\ldots,Y_n$ be an
i.i.d. sample of real random variables with distribution $G$ (not
necessarily continuous), and
let $\hat{G}_n$ be the corresponding empirical distribution function.
Then
\[\mbox{P}\left(\sup_{t\in\re}\sqrt{n}\left|\hat{G}_n(t)-G(t)\right|
>x\right)\leq 2e^{-2x^2},\]
for all $x\geq 0$.
\end{theorem}

Proof. This is the celebrated result of Dvoretsky, Kiefer and
Wolfowitz (1956), given in their lemma~2, as
refined by Massart (1990) in his corollary~1. 
We omit the proof of their result
but note that their result applies to the special
case where $G$ is continuous. We now show that it also applies
when $G$ may be discontinuous. Without loss of generality, assume
that $G$ has discontinuities, and let $T_1,\ldots,T_m$ be the 
locations of the discontinuities of $G$, where $m$ may be infinity. 
Note that the number of discontinuities
can be at most countable. Let $r_1,\ldots,r_m$ be the jump sizes
of $G$ at $T_1,\ldots,T_m$. Now let $U_1,\ldots,U_n$ be i.i.d.
uniform random variables independent of the $Y_1,\ldots,Y_n$,
and define new random variables 
$Z_i=Y_i+\sum_{j=1}^mr_j\left[\ind\{T_j<Y_i\}+\ind\{T_j=Y_i\}U_i\right]$,
$1\leq i\leq n$. Define also the transformation 
$t\mapsto R(t)=t+\sum_{j=1}^mr_j\ind\{T_j\leq t\}$; let
$\hat{H}_n^{\ast}$ be the empirical distribution of $Z_1,\ldots,Z_n$; and 
let $H$ be the distribution of $Z_1$. It is not hard to verify that
\begin{eqnarray*}
\sup_{t\in\re}|\hat{G}_n(t)-G(t)|&=&\sup_{t\in\re}|\hat{H}_n(R(t))
-H(R(t))|\\
&\leq&\sup_{s\in\re}|\hat{H}_n(s)-H(s)|,
\end{eqnarray*}
and the desired result now follows since $H$ is continuous.$\Box$ 

{\it Proof of theorem~\ref{t1}.} 
Let $U_{ij(n)}$, $i=0,\ldots,n$ and
$j=1,\ldots,p_n$, be independent uniform random variables.
Then, by theorem~\ref{t2} below, there exist Brownian bridges
$B_{1(n)},\ldots,B_{p_n(n)}$, where, for each $1\leq j\leq p_n$,
$B_{j(n)}$ depends only on $X_{1j(n)},\ldots,X_{nj(n)}$ and
$U_{0j(n)},\ldots,U_{nj(n)}$ and
\begin{eqnarray}
P\left(\sqrt{n}\left\|\sqrt{n}(\hat{F}_{j(n)}-F_{j(n)})-B_{j(n)}(F_{j(n)})
\right\|_{\infty}
>x+12\log n\right)&\leq& 2 e^{-x/6},\label{t1.e2}
\end{eqnarray}
for all $x\geq 0$ and all $n\geq 2$. Now define
\[U_{j(n)}=\left(\frac{\sqrt{n}}{\log n}\left\|\sqrt{n}(\hat{F}_{j(n)}
-F_{j(n)})-B_{j(n)}(F_{j(n)})\right\|_{\infty}-12\right)^{+},\]
where $u^{+}$ is the positive part of $u$. By lemma~2.2.1 of VW,
expression~(\ref{t1.e2}) implies
that $\|U_{j(n)}\|_{\psi_1}\leq 18/\log n$. Reapplying the
result that $\log(k+1)\leq 2\log k$
for any $k\geq 2$, we now have, by the fact that
$\lim\sup_{x,y\rightarrow\infty}\psi_1(x)\psi_1(y)/\psi_1(xy)=0$
combined with lemma~2.2.2 of VW, that there exists a universal
constant $0<c_2<\infty$ for which
\[\left\|\max_{1\leq j\leq p_n}U_{j(n)}\right\|_{\psi_1}
\leq \frac{c_2\log p_n}{\log n}.\]
Now~(\ref{t1.e1}) follows, for $c_1=12$,
from the definition of $U_{j(n)}$.$\Box$
 
\begin{theorem}\label{t2}
For $n\geq 2$,
let $Y_1,\ldots,Y_n$ be i.i.d. real random variables with
distribution $G$ (not necessarily continuous),
and let $U_0,\ldots,U_n$ be independent uniform
random variables independent of $Y_1,\ldots,Y_n$. Then there
exists a standard Brownian motion $B$ depending only on
$Y_1,\ldots,Y_n$ and $U_0,\ldots,U_n$ such that, for all $x\geq 0$,
\begin{eqnarray}
P\left(\sqrt{n}\left\|\sqrt{n}(\hat{G}_n-G)-B(G)\right\|_{\infty}
>x+12\log n\right)&\leq& 2 e^{-x/6},\label{t2.e1}
\end{eqnarray}
where $\hat{G}_n$ is the empirical distribution of $Y_1,\ldots,Y_n$.
\end{theorem}

Proof. We will apply the same method for handling the discontinuities of 
$G$ as used in the proof of theorem~\ref{t.c1.t1}.
Let $m\geq 0$, $T_1,\ldots,T_m$, and
$r_1,\ldots,r_m$ be as defined in the proof of theorem~\ref{t.c1.t1}. 
Similarly define $Z_1,\ldots,Z_m$, $R$, $\hat{H}_n$ and $H$,
except that we will utilize the uniform random variables 
$U_1,\ldots,U_n$ given in the statement of theorem~\ref{t2}. 
By the continuity of $H$ as established in the proof of 
theorem~\ref{t.c1.t1}, $H(Z_1)$ is now uniformly distributed.
Thus, by the Hungarian construction theorem (theorem~1) 
of Bretagnolle and Massart (1989), there exists
a Brownian bridge $B$ depending only on $Z_1,\ldots,Z_n$ and $U_0$
such that
\[P\left(\sqrt{n}\left\|\sqrt{n}(\hat{H}_n-H)-B(H)\right\|_{\infty}
>x+12\log n\right)\leq 2 e^{-x/6},\]
for all $x\geq 0$. The desired result now follows since
\begin{eqnarray*}
\sup_{t\in\re}\left|\sqrt{n}(\hat{G}_n(t)-G(t))-B(G(t))\right|
&=&\sup_{t\in\re}\left|\sqrt{n}(\hat{H}_n(R(t))-H(R(t)))
-B(H(R(t)))\right|\\
&\leq&\sup_{s\in\re}\left|\sqrt{n}(\hat{H}_n(s)-H(s))-B(H(s))
\right|.\Box
\end{eqnarray*}

\vspace{0.35cm}
{\it Proof of corollary~\ref{c3}}.
The result is a consequence of theorem~\ref{t.c1} via the
following integration by parts identity:
\begin{eqnarray}
&&\label{c3.e2}\\
\int_{[a_{j(n)},b_{j(n)}]}x\left[
d\hat{F}_{j(n)}(x)-dF_{j(n)}(x)\right]&=&-\int_{[a_{j(n)},b_{j(n)}]}
\left[\hat{F}_{j(n)}(x)-F_{j(n)}(x)\right]dx.\Box\nonumber
\end{eqnarray}

\vspace{0.35cm}
{\it Proof of corollary~\ref{c2}}.
Note that for any $x\in\re$ and any $y>0$,
\begin{eqnarray*}
|\Phi(xy)-\Phi(x)|&\leq&\sup_{1\wedge y\leq u\leq 1\vee y}|x|\phi(xu)|y-1|\\
&\leq&0.25\times\sup_{1\wedge y\leq u\leq 1\vee y}\frac{|y-1|}{u}\\
&\leq&0.25\times|1-y|\vee\left|1-\frac{1}{y}\right|.
\end{eqnarray*}
The constant 0.25 comes from the fact that
$\sup_{u>0}u\phi(u)\leq(2\pi e)^{-1/2}\leq 0.25$. Thus
\begin{eqnarray}
\max_{1\leq j\leq p_n}
|\hat{\pi}_{j(n)}-\hat{\pi}_{j(n)}^{\ast}|
&\leq&\frac{1}{2}
\left(
\max_{1\leq j\leq n}
(\hat{\sigma}_{j(n)}\vee\sigma_{j(n)})
\left|\frac{1}{\hat{\sigma}_{j(n)}}-\frac{1}
{\sigma_{j(n)}}\right|\right),\label{c2.e4}
\end{eqnarray}
where $\hat{\pi}_{j(n)}^{\ast}=2\Phi(-|T_{j(n)}^{\ast}|)$ and
$T_{j(n)}^{\ast}=\sqrt{n}(\bar{X}_{j(n)}-\mu_{0,j(n)})/
\sigma_{j(n)}$.

Now the integration by parts formula~(\ref{c3.e2}) combined with
theorem~\ref{t1} yields
\begin{eqnarray*}
\lefteqn{\left\|\max_{1\leq j\leq p_n}\left|\frac{\sqrt{n}(\bar{X}_{j(n)}
-\mu_{j(n)})}{\sigma_{j(n)}}+\frac{\int_{[a_{j(n)},b_{j(n)}]}
B_{j(n)}(F_{j(n)}(x))dx}{\sigma_{j(n)}}\right|\;\right\|_{\psi_1}}&&\\
&\mbox{\hspace{2.0in}}&\leq\left(\frac{c_1\log n+c_2\log p_n}{\sqrt{n}}\right)
\max_{1\leq j\leq p_n}\frac{|b_{j(n)}-a_{j(n)}|}{\sigma_{j(n)}},
\end{eqnarray*} 
where $c_1,c_2$ and $B_{1(n)},\ldots,B_{p_n(n)}$ are as given in
theorem~\ref{t1}, and where
$$Z_{j(n)}=-\int_{[a_{j(n)},b_{j(n)}]}B_{j(n)}(F_{j(n)}(x))dx
/\sigma_{j(n)}$$ is standard normal for all $1\leq j\leq p_n$.
This, combined with the fact
that $|\Phi(x)-\Phi(y)|\leq |x-y|/2$ for all $x,y\in\re$, yields
the desired result.$\Box$

\vspace{0.35cm}
{\it Proof of corollary~\ref{c4}}.
That the left-hand-side of~(\ref{c4.e2}) is $o_P(1)$ follows
from condition~(\ref{c4.e1}) combined with theorem~\ref{t.c1}.
By the definition of the sample median, we have that
$\hat{F}_{j(n)}(\hat{\xi}_{j(n)})-F_{j(n)}(\xi_{j(n)})\equiv E_{j(n)}$, where
$|E_{j(n)}|\leq 1/n$. This now implies that
$\hat{F}_{j(n)}(\hat{\xi}_{j(n)})-F(\hat{\xi}_{j(n)})+
F(\hat{\xi}_{j(n)})-F(\xi_{j(n)})=E_{j(n)}$. The result now follows
from the mean value theorem and condition~(\ref{c4.e1}).$\Box$

\vspace{0.35cm}
{\it Proof of Corollary \ref{c5}}.
Now, for some
$\xi_{j(n)}^{\ast}$ in between $\xi_{j(n)}$ and $\hat{\xi}_{j(n)}$,
we have $f_{j(n)}(\xi_{j(n)}^{\ast})(\hat{\xi}_{j(n)}-\xi_{j(n)})
=-F_{j(n)}(\hat{\xi}_{j(n)})+F_{j(n)}(\xi_{j(n)})$.
Using the conditions of the corollary, we obtain that the
$\hat{f}_{j(n)}$ terms are simultaneously consistent for
the quantities $f_{j(n)}(\xi_{j(n)}^{\ast})$ and that these later
quantities are bounded above and below. Now we can argue as
in the first part of the proof of corollary~\ref{c2} that
$\max_{1\leq j\leq p_n}|\hat{\pi}_{j(n)}'- 
\tilde{\pi}_{j(n)}'|=o_P(1)$, where $\tilde{\pi}_{j(n)}'=
2\Phi(-|\tilde{U}_{j(n)}|)$ and
\begin{eqnarray*}
\tilde{U}_{j(n)}&=&2\sqrt{n}f_{j(n)}(\xi_{j(n)}^{\ast})
(\hat{\xi}_{j(n)}-\xi_{0,j(n)})\\
&=&-2\sqrt{n}(F_{j(n)}(\hat{\xi}_{j(n)}) 
-F_{j(n)}(\xi_{j(n)}))+2\sqrt{n}f_{j(n)}(\xi_{j(n)}^{\ast})
(\xi_{j(n)}-\xi_{0,j(n)}).
\end{eqnarray*}

Note that 
\begin{eqnarray*}
\lefteqn{\sqrt{n}\left(F_{j(n)}(\hat{\xi}_{j(n)})-
F_{j(n)}(\xi_{j(n)})\right)}\mbox{\hspace{0.5in}}&&\\
&=&-\sqrt{n}\left(\hat{F}_{j(n)}(\hat{\xi}_{j(n)})
-F_{j(n)}(\hat{\xi}_{j(n)})-\hat{F}_{j(n)}(\xi_{j(n)})
+F_{j(n)}(\xi_{j(n)})\right)\\
&&-\sqrt{n}\left(\hat{F}_{j(n)}(\xi_{j(n)})-F_{j(n)}(\xi_{j(n)})\right)
+\sqrt{n}\left(\hat{F}_{j(n)}
(\hat{\xi}_{j(n)})-F_{j(n)}(\xi_{j(n)})\right)\\
&\equiv&-A_{j(n)}-V_{j(n)}+C_{j(n)},
\end{eqnarray*}
where $C_{j(n)}=\sqrt{n}E_{j(n)}$ and $E_{j(n)}$ as defined in the 
proof of corollary~\ref{c4} with $|E_{j(n)}|\leq 1/n$.
Hence $C_{j(n)}$ vanishes asymptotically, uniformly over $1\leq j\leq p_n$.
Theorem~\ref{t1} tells us that we can, uniformly over $1\leq j\leq p_n$,
replace $A_{j(n)}$ and $V_{j(n)}$ with
$A_{j(n)}'=B_{j(n)}(F_{j(n)}(\hat{\xi}_{j(n)}))-
B_{j(n)}(F_{j(n)}(\xi_{j(n)}))$ and $V_{j(n)}'=B_{j(n)}(1/2)$.
Note that $Z_{j(n)}\equiv 2B_{j(n)}(1/2)$ are standard normals and
and that $B_{j(n)}(t)=W_{j(n)}(t)-tW_{j(n)}(1)$, for all $t\in[0,1]$,
for some standard Brownian motions $W_{j(n)}$. Thus, by the
symmetry properties of Brownian motion, $|A_{j(n)}'|$
\[\leq\sqrt{\hat{\delta}_{j(n)}}\left[
\sup_{0\leq t\leq\hat{\delta}_{j(n)}}|W_{j(n)}'(t)|
+\sup_{0\leq t\leq\hat{\delta}_{j(n)}}|W_{j(n)}''(t)|\right]
+\hat{\delta}_{j(n)}|W_{j(n)}(1)|\equiv
\tilde{A}_{j(n)}(\hat{\delta}_{j(n)}),\] 
where $\hat{\delta}_{j(n)}\equiv M|\hat{\xi}_{j(n)}-\xi_{j(n)}|$;
$M$ is as defined in~(\ref{c4.e2}); and where $W_{j(n)}$, $W_{j(n)}'$
and $W_{j(n)}''$ are Brownian motions.

Now, for each $k<\infty$ and $\rho>0$, we have
\begin{eqnarray}
\lefteqn{\pr\left(\max_{1\leq j\leq p_n}|A_{j(n)}'|>\rho\right)}
\mbox{\hspace{0.5in}}&&\label{c4.p.e1}\\
&\leq&
\pr\left(\max_{1\leq j\leq p_n}\tilde{A}_{j(n)}(kr_n)>\rho\right)
+\pr\left(\max_{1\leq j\leq p_n}\hat{\delta}_{j(n)}>kr_n\right),\nonumber
\end{eqnarray}
where $r_n\equiv\log(n\vee p_n)/n+\sqrt{\log p_n /n}$. However,
using the facts that a standard normal deviate and the 
supremum of the absolute value of a Brownian motion over $[0,1]$
both have sub-Gaussian tails (i.e., have bounded $\psi_2$-norms), 
we have $\max_{1\leq j\leq p_n}\tilde{A}_{j(n)}(kr_n)
\leq O_P\left(\sqrt{\log p_n}\left[r_n+\sqrt{r_n}\right]\right)
\rightarrow 0$, in probability, since $\log^3 p_n/n\rightarrow 0$.
Thus the first term on the right-hand-side of~(\ref{c4.p.e1})
goes to zero. Since corollary~\ref{c4} implies
$\lim_{k\rightarrow\infty}\limsup_{n\rightarrow\infty}
\pr\left(\max_{1\leq j\leq p_n}\hat{\delta}_{j(n)}>kr_n\right)=0$,
the left-hand-side of~(\ref{c4.p.e1}) also goes to zero
as $n\rightarrow\infty$. Thus $\tilde{U}_{j(n)}$ can be approximated by
$U_{j(n)}'=Z_{j(n)}+2\sqrt{n}f_{j(n)}(\xi_{j(n)}^{\ast})
(\xi_{j(n)}-\xi_{0,j(n)})$
simultaneously over all $1\leq j\leq p_n$.

Now we can use arguments given at the beginning of the proof of
corollary~\ref{c2} (again) in combination with the simultaneous
consistency of $\hat{\xi}_{j(n)}$ and the assumed properties
of $f_{j(n)}$ to obtain that
\[\max_{1\leq j\leq p_n}\left|\Phi(-|U_{j(n)}'|)
-\Phi\left(-\left|\frac{f_{j(n)}(\xi_{j(n)})}{f_{j(n)}(\xi_{j(n)}^{\ast})}
Z_{j(n)}+2\sqrt{n}f_{j(n)}(\xi_{j(n)})(\xi_{j(n)}-\xi_{0,j(n)})\right|
\right)\right|\]
$=o_P(1)$. Now define $\eta_n\equiv\max_{1\leq j\leq p_n}|\hat{\xi}_{j(n)}
-\xi_{j(n)}|$. By condition~(\ref{c5.e2}), we have that
\begin{eqnarray*}
\lefteqn{\max_{1\leq j\leq p_n}\left|\frac{f_{j(n)}(\xi_{j(n)})}
{f_{j(n)}(\xi_{j(n)}^{\ast})}-1\right|\times|Z_{j(n)}|\;\;\leq\;\;
O_P\left(\max_{1\leq j\leq p_n}|Z_{j(n)}|\right)}\mbox{\hspace{1.0in}}&&\\
&&\times
\max_{1\leq j\leq p_n}\sup_{\epsilon\leq\eta_n}\sup_{u:|u|\leq\epsilon}
\frac{|f_{j(n)}(\xi_{j(n)}+u)-f_{j(n)}(\xi_{j(n)})|}{\epsilon^{1/2}}
\eta_n^{1/2}\\
&\leq&O_P(\sqrt{\log p_n})\times\eta_n^{1/2}\\
&=& O_P\left(\sqrt{\log p_n}\times\sqrt{\frac{\log n\vee p_n}{n}
+\sqrt{\frac{\log p_n}{n}}}\right)\\
&\rightarrow&0,
\end{eqnarray*}
where the equality follows from corollary~\ref{c4}.
The desired result now follows.$\Box$

{\it Proof of corollary~\ref{c6}.} The proof follows the same general
logic as the proof of corollary~\ref{c2}. Using the fact that,
for any $x\in\re$ and any $y>0$,
$|\Phi(xy)-\Phi(x)|\leq 0.25\times |1-y|\vee |1-y^{-1}|$, we have
\begin{eqnarray}
\;\;\max_{1\leq j\leq p_n}|\hat{\pi}_{j(n)}^{\ast}
-\hat{\pi}_{j(n)}^{\ast\ast}|
&\leq&\max_{1\leq j\leq p_n}\frac{1}{2}\left\{\left|
\left(\frac{n_1\left[\sigma_{j(n)}^{(2)}\right]^2
+n_2\left[\sigma_{j(n)}^{(1)}\right]^2}{n_1\left[\hat{\sigma}_{j(n)}^{(2)}
\right]^2+n_2\left[\hat{\sigma}_{j(n)}^{(1)}\right]^2}\right)^{1/2}-1\right|
\right.\label{c6.pe1}\\
&&\vee\left.\left|\left(\frac{n_1\left[\hat{\sigma}_{j(n)}^{(2)}
\right]^2+n_2\left[\hat{\sigma}_{j(n)}^{(1)}\right]^2}
{n_1\left[\sigma_{j(n)}^{(2)}\right]^2+n_2\left[\sigma_{j(n)}^{(1)}\right]^2}
\right)^{1/2}-1\right|\right\},\nonumber
\end{eqnarray}
where $\hat{\pi}_{j(n)}^{\ast\ast}=2\Phi(-|T_{j(n)}^{\ast\ast}|)$ and
\begin{eqnarray*}
T_{j(n)}^{\ast\ast}&=&\left(\frac{n_2\left[\sigma_{j(n)}^{(1)}\right]^2}
{n_1\left[\sigma_{j(n)}^{(2)}\right]^2+n_2\left[\sigma_{j(n)}^{(1)}\right]^2}
\right)^{1/2}\frac{\sqrt{n_1}\left(\bar{X}_{j(n)}^{(1)}-\mu_{j(n)}^{(1)}
\right)}{\sigma_{j(n)}^{(1)}}\\
&&-\left(\frac{n_1\left[\sigma_{j(n)}^{(2)}
\right]^2}{n_1\left[\sigma_{j(n)}^{(2)}\right]^2+n_2\left[\sigma_{j(n)}^{(1)}
\right]^2}\right)^{1/2}\frac{\sqrt{n_2}\left(\bar{X}_{j(n)}^{(2)}
-\mu_{j(n)}^{(2)}\right)}{\sigma_{j(n)}^{(2)}}\\
&&+\sqrt{\frac{n_1 n_2}{n_1\left[\sigma_{j(n)}^{(2)}\right]^2
+n_2\left[\sigma_{j(n)}^{(1)}\right]^2}}
\left(\mu_{j(n)}^{(1)}-\mu_{j(n)}^{(2)}\right).
\end{eqnarray*}
Now, virtually identical Brownian bridge approximation arguments to those 
used in the proof of corollary~\ref{c2} yield that
\[\max_{1\leq j\leq p_n}\left|\hat{\pi}_{j(n)}^{\ast\ast}-\pi_{j(n)}^{\ast}
\right|\leq\sum_{k=1,2}\frac{c_1\log n_k+c_2\log p_n}
{\sqrt{n_k}}\left(\max_{1\leq j\leq p_n}\frac{|b_{j(n)}^{(k)}-a_{j(n)}^{(k)}|}
{\sigma_{j(n)}^{(k)}}\right).\]

In order to finish the proof, we need to bound the right-hand-side
of~(\ref{c6.pe1}). To begin with, note that for any scalars
$c_1,c_2,d_1,d_2\geq 0$,
\begin{eqnarray*}
\left|\left(\frac{n_1c_2^2+n_2c_1^2}{n_1d_2^2+n_2d_1}
\right)^{1/2}-1\right|
&\leq&\left|\left(\frac{n_1c_2^2+n_2c_1^2}{n_1d_2^2+n_2d_1}\right)^{1/2}
-\left(\frac{n_1c_2^2+n_2d_1^2}{n_1d_2^2+n_2d_1}\right)^{1/2}\right|\\
&&+\left|\left(\frac{n_1c_2^2+n_2d_1^2}{n_1d_2^2+n_2d_1}\right)^{1/2}
-\left(\frac{n_1d_2^2+n_2d_1^2}{n_1d_2^2+n_2d_1}\right)^{1/2}\right|\\
&\leq&\left(\frac{n_2d_1^2}{n_1d_2^2+n_2d_1^2}\right)^{1/2}
\left|\frac{c_1}{d_1}-1\right|\\
&&+\left(\frac{n_1d_2^2}{n_1d_2^2+n_2d_1^2}\right)^{1/2}
\left|\frac{c_2}{d_2}-1\right|\\
&\leq&\left|\frac{c_1}{d_1}-1\right|+\left|\frac{c_2}{d_2}-1\right|,
\end{eqnarray*}
where the second inequality follows from the fact that for any
$a,b,x,y\geq 0$, 
\[\left|\left(ax^2+b\right)^{1/2}-\left(ay^2+b\right)^{1/2}
\right|\leq\sqrt{a}|x-y|.\] 
Hence both
\[\left|\left(\frac{n_1\left[\sigma_{j(n)}^{(2)}\right]^2
+n_2\left[\sigma_{j(n)}^{(1)}\right]^2}{n_1\left[\hat{\sigma}_{j(n)}^{(2)}
\right]^2+n_2\left[\hat{\sigma}_{j(n)}^{(1)}\right]^2}\right)^{1/2}-1\right|
\;\leq\;\left|\frac{\sigma_{j(n)}^{(1)}}{\hat{\sigma}_{j(n)}^{(1)}}-1\right|
+\left|\frac{\sigma_{j(n)}^{(2)}}{\hat{\sigma}_{j(n)}^{(2)}}-1\right|\]
and
\[\left|\left(\frac{n_1\left[\hat{\sigma}_{j(n)}^{(2)}\right]^2
+n_2\left[\hat{\sigma}_{j(n)}^{(1)}\right]^2}{n_1\left[\sigma_{j(n)}^{(2)}
\right]^2+n_2\left[\sigma_{j(n)}^{(1)}\right]^2}\right)^{1/2}-1\right|
\;\leq\;\left|\frac{\hat{\sigma}_{j(n)}^{(1)}}
{\sigma_{j(n)}^{(1)}}-1\right|+\left|\frac{\hat{\sigma}_{j(n)}^{(2)}}
{\sigma_{j(n)}^{(2)}}-1\right|,\]
and thus the right-hand-side of~(\ref{c6.pe1}) is bounded by
\[\frac{1}{2}\sum_{k=1,2}\max_{1\leq j\leq p_n}
\left(\hat{\sigma}_{j(n)}^{(k)}\vee\sigma_{j(n)}^{(k)}\right)
\left|\frac{1}{\hat{\sigma}_{j(n)}^{(k)}}-\frac{1}
{\sigma_{j(n)}^{(k)}}\right|,\]
completing the proof.$\Box$

{\it Proof of corollary~\ref{c7}.} The proof consists of extending
the proof of corollary~\ref{c5} in a manner similar to the way 
in which the proof of corollary~\ref{c2} was extended for proving 
corollary~\ref{c6}. A key difference is that the role of $\sigma_{j(n)}^{(k)}$
and $\hat{\sigma}_{j(n)}^{(k)}$ in the proof of corollary~\ref{c6}
is replaced by $1/f_{j(n)}^{(k)}$ and $1/\hat{f}_{j(n)}^{(k)}$,
for $k=1,2$ and $1\leq j\leq p_n$. The remaining necessary
extensions of the proof of corollary~\ref{c5} are straightforward.$\Box$

{\it Proof of lemma~\ref{l2}.} Because of the invariance under monotone
transformation, we can assume without loss of generality that the
data are uniformly distributed. Classical arguments in Billingsley
(1968) yield the second result. In particular, the form of
the limiting distribution function, which is the distribution
of the supremum in absolute value of a Brownian bridge, 
can be found on page~85 of Billingsley.
Arguments for establishing the remaining two results can be
found in section~3.9.4 (for the Wilcoxon statistic) and in
section~2.13.2 (for the Cram\'{e}r-von Mises statistic) of
van der Vaart and Wellner (1996).$\Box$

{\it Proof of theorem~\ref{t3}.} Define $\tilde{E}_n=
\max_{1\leq j\leq p_n}\left\|\tilde{F}_{j(n)}-\hat{F}_{j(n)}\right\|_{\infty}$
and, for each $\delta\geq 0$, $\hat{E}_n(\delta)=
\max_{1\leq j\leq p_n}\sup_{|s-t|\leq\delta}\left|\tilde{F}_{j(n)}(s)
-\hat{F}_{j(n)}(t)\right|$. Suppose now that for some positive,
non-increasing sequences $\{s_n,\delta_n\}$, with $\delta_n\rightarrow 0$, 
we have $\hat{E}_n(\delta_n)=o_P(s_n)$ and 
$\pr(\hat{\epsilon}_n>\delta_n)=o(1)$. Then, by the definition of
$\hat{\epsilon}_n$, 
\begin{eqnarray}
\tilde{E}_n&=&\tilde{E}_n\ind\{\hat{\epsilon}_n\leq \delta_n\}
+\tilde{E}_n\ind\{\hat{\epsilon}_n>\delta_n\}
\;\;\leq\;\;\hat{E}_n(\delta_n)+o_P(s_n)
\;\;=\;\;o_P(s_n).\label{t3.e1}
\end{eqnarray}
Now, by theorem~\ref{t1} and condition~(\ref{t3.e0}),  
we have for any sequence $\delta_n\downarrow 0$,
\begin{eqnarray*}
\sqrt{n}\hat{E}_n(\delta_n)&\leq&\max_{1\leq j\leq p_n}
\sup_{|s-t|\leq\delta_n}\sqrt{n}\left|
\hat{F}_{j(n)}(s)-F_{j(n)}(s)-\hat{F}_{j(n)}(t)+F_{j(n)}(t)
\right|+\sqrt{n}\tilde{M}\delta_n\\
&\leq&\max_{1\leq j\leq p_n}\sup_{|s-t|\leq\delta_n}
\left|B_{j(n)}(F_{j(n)}(s))-B_{j(n)}(F_{j(n)}(t))\right|\\
&&+O_P\left(\frac{\log n+\log p_n}{\sqrt{n}}+\sqrt{n}\delta_n\right).
\end{eqnarray*}
Combining this with a reapplication of condition~(\ref{t3.e0}) along
with lemma~\ref{l1} below (a precise modulus of continuity bound
for Brownian motion), we obtain
\begin{eqnarray}
\sqrt{n}\hat{E}_n(\delta_n)&\leq&
O_P\left(\sqrt{(\log p_n)\delta_n\log(1/\delta_n)}
+\frac{\log n+\log p_n}{\sqrt{n}}+\sqrt{n}\delta_n\right).\label{t3.e2}
\end{eqnarray}
Both~(\ref{t3.e1}) and~(\ref{t3.e2}) will prove useful at several points 
in our proof.

Using the fact that $\hat{\epsilon}_n=o_P(1)$, we can find a positive,
sufficiently slowly decreasing sequence $\delta_n\rightarrow 0$ such that 
$\hat{\epsilon}_n=o_P(\delta_n)$. Now, by applying~(\ref{t3.e1}) with
$s_n=1$, we obtain result~(i) of the 
theorem: $\tilde{E}_n=o_P(1)$. For result~(ii), we can use the
fact that $\log^2 p_n/n=o(1)$, to construct a positive,
non-decreasing sequence $r_n\rightarrow\infty$ slowly enough so that
$r_n\log p_n/\sqrt{n}=o(1)$ and $r_n/\log n=o(1)$. Since
$\sqrt{n}(\log n)\hat{\epsilon}_n=O_P(1)$, we have
\[\sqrt{n}\hat{\epsilon}_n\log(1/\hat{\epsilon}_n)
=\sqrt{n}(\log n)\hat{\epsilon}_n\left(\frac{\log\sqrt{n}
-\log(\sqrt{n}\hat{\epsilon}_n)}{\log n}\right)=O_P(1).\]
Thus, if we set $\delta_n=r_n/(\sqrt{n}\log n)$, we have
$\hat{\epsilon}_n=o_P(\delta_n)$. We also have, by~(\ref{t3.e2}), that
\[\hat{E}_n(\delta_n)=O_P\left(\sqrt{\frac{1}{n}\times
\frac{r_n\log p_n}{\sqrt{n}}\times\frac{\log\sqrt{n}+\log\log n-\log r_n}
{\log n}}+o\left(\frac{1}{\sqrt{n}}\right)\right)
=o_P(n^{-1/2}).\]
The proof is now complete by reapplying~(\ref{t3.e1})
with the choice $s_n=n^{-1/2}$.$\Box$

\begin{lemma}\label{l1}
Let $W:[0,1]\mapsto\re$ be a standard Brownian motion. Then there
exists a universal constant $k_0<\infty$ such that
\[\left\|\sup_{|s-t|\leq\delta}|W(s)-W(t)|\,\right\|_{\psi_2}
\leq k_0\sqrt{\delta\log(1/\delta)}\]
for all $0<\delta\leq 1/2$.
\end{lemma}

{\it Proof.} Fix $\delta\in(0,1)$. Let $n_{\delta}$ be the smallest
integer $\geq 1+1/\delta$, and extend the Brownian
motion $W$ to the interval $[0,\delta n_{\delta}]$. Now
\begin{eqnarray}
\sup_{|s-t|\leq\delta}|W(s)-W(t)|&\leq&\max_{1\leq j\leq n_{\delta}}
\sup_{(j-1)\delta\leq s<t\leq(j+1)\delta}|W(s)-W(t)|\label{l1.e1}\\
&\leq&2\max_{1\leq j\leq n_{\delta}}\sup_{t\in[(j-1)\delta,(j+1)\delta]}
|W(t)-W((j-1)\delta)|\nonumber\\
&\leq&2\max_{1\leq j\leq n_{\delta}}\sup_{t\in[0,1]}\sqrt{2\delta}
|W_j^{\ast}(t)|,\nonumber
\end{eqnarray}
where $W_1^{\ast},\ldots,W_{n_{\delta}}^{\ast}$ are a dependent
collection of standard Brownian motions. The last inequality
follows from the symmetry properties of Brownian motion. We can now use the
fact that the tail probabilities of the supremum over $[0,1]$ of
the absolute value of Brownian motion are sub-Gaussian (and 
thus have bounded $\psi_2$-norms) to obtain that the 
$\psi_2$-norm of the left side of~(\ref{l1.e1}) is bounded by
$k_{\ast}2\sqrt{2\delta\log(1+n_{\delta})}\leq k_{\ast}2\sqrt{2\delta
\log(3+1/\delta)}\leq k_0\sqrt{\delta\log(1/\delta)}$,
where $k_0=5k_{\ast}$ does not depend on $\delta$. The last
inequality follows because
$\log(3+1/\delta)/\log(1/\delta)\leq 1+\log(1+3\delta)/\log(1/\delta)\leq 3$ 
for all $\delta\in(0,1/2]$.$\Box$

{\em Proof of corollary~\ref{c8}.} Result~(i) follows directly from
part~(i) of theorem~\ref{t3} and theorem~\ref{t.c1}. Result~(ii) is a direct
consequence of part~(i) of theorem~\ref{t3} and a minor modification of the 
integration by parts identity~(\ref{c3.e2}) used in the proof of 
corollary~\ref{c2}. The proof of result~(iii) is a straightforward
extension of the proof of corollary~\ref{c4} which incorporates the
conclusion of part~(i) of theorem~\ref{t3}.$\Box$

{\em Proof of corollary~\ref{c9}.} For result~(i), we use part~(ii)
of theorem~\ref{t3} combined with integration
by parts to obtain that 
\[\max_{1\leq j\leq p_n}\left|\check{X}_{j(n)}
-\bar{X}_{j(n)}\right|=o_P\left(n^{-1/2}\left[\max_{1\leq j\leq p_n}|b_{j(n)}
-a_{j(n)}|+2\hat{\epsilon}_n\right]\right)=o_P(1).\] 
Now corollary~\ref{c3} gives us the desired results since 
\[\sqrt{\frac{\log p_n}{n}}\max_{1\leq j\leq p_n}|b_{j(n)}-a_{j(n)}|
=\sqrt{\frac{\log p_n}{\sqrt{n}}}\max_{1\leq j\leq p_n}n^{1/4}
|b_{j(n)}-a_{j(n)}|=o(1).\]
For result~(ii), we also use part~(ii) of theorem~\ref{t3} combined
with integration by parts to obtain
\[\max_{1\leq j\leq p_n}\left|\frac{\sqrt{n}\left(
\check{X}_{j(n)}-\bar{X}_{j(n)}\right)}{\sigma_{j(n)}}\right|
\leq o_P\left(\max_{1\leq j\leq p_n}\frac{|b_{j(n)}-a_{j(n)}|+
2\hat{\epsilon}_n}{\sigma_{j(n)}}\right)=o_P(1),\]
and the desired result follows using the Brownian bridge approximation
of $\sqrt{n}\left(\bar{X}_{j(n)}\right.$ 
$\left.-\mu_{j(n)}\right)/\sigma_{j(n)}$
given in the proof of corollary~\ref{c2}. 
For result~(iii), the desired conclusion is obtained via part~(ii) of 
theorem~\ref{t3} combined with
a straightforward adaptation of the proof of corollary~\ref{c5}.$\Box$

{\em Proof of corollary~\ref{c10}.} The proof follows almost immediately
from applying part~(ii) of theorem~\ref{t3} to each sample separately,
yielding the result 
\[\max_{k=1,2}\,\max_{1\leq j\leq p_n}\sqrt{n_k}\left\|
\tilde{F}_{j(n)}^{(k)}-\hat{F}_{j(n)}^{(k)}\right\|_{\infty}=o_P(1).\]
Now, the proofs of results~(i) and~(ii)
are direct extensions of the one-sample results of corollary~\ref{c9}
combined with straightforward adaptations of arguments found in the
proofs of corollaries~\ref{c6} and~\ref{c7}. The proof of result~(iii)
also follows almost immediately. For the Kolmogorov-Smirnov statistic,
the result is obvious. For the other two statistics, the result follows
with some help from integration by parts.$\Box$
\vspace{0.1in}

\begin{center}
REFERENCES
\end{center}
\begin{list}{}{\setlength{\leftmargin}{0.2in}
\setlength{\itemindent}{-.2in}}

\item{\sc Benjamini, Y.}, and {\sc Hochberg, Y.} (1995). Controlling
the false discovery rate: A practical and powerful approach to
multiple testing. {\em Journal of the Royal Statistical Society},   
Series B, {\bf 57}, 289--300.

\item{\sc Billingsley, P.} (1968). {\em Convergence of Probability
Measures}. Wiley, New York.

\item{\sc Bretagnolle, J.}, and {\sc Massart, P.} (1989). Hungarian
construction from the nonasymptotic viewpoint. {\em Annals of
Probability} {\bf 17}, 239--256.

\item{\sc Dudoit, S., Fridlyand, J.} and {\sc Speed, T. P.} (2002).
Comparison of discrimination methods for the classification of tumors
using gene expression data. {\em Journal of the American Statistical
Association} {\bf 97}, 77--87.

\item{\sc Dvoretsky, A., Kiefer, J.}, and {\sc Wolfowitz, J.} (1956).
Asymptotic minimax character of the sample distribution function and
of the classical multinomial estimator. {\em Annals of Mathematical
Statistics} {\bf 27}, 642--669.

\item{\sc Fan, J., Peng, H.} and {\sc Huang, T.} (2004) Semilinear
high-dimensional model for normalization of microarray data: a
theoretical analysis and partial consistency. {\it Journal of the American 
Statistical Association}, In press.

\item{\sc Ghosh, D.} and {\sc Chinnaiyan, A.M.} (2004)
Classification and selection of biomarkers in genomic data using
LASSO. {\it Journal of Biomedicine and Biotechnology}, In press.

\item{\sc Gui,  J.} and {\sc Li,  H.} (2004).
Penalized Cox Regression Analysis in the High-Dimensional
and Low-sample Size Settings, with Applications to Microarray 
Gene Expression Data. Submitted.

\item{\sc Huang, J., Wang, D.}, and {\sc Zhang, C.-H.} (2005).
A two-way semi-linear model for normalization and analysis of 
cDNA microarray data. {\em Journal of the American Statistical 
Association}, In press.

\item{\sc Koml\'{o}s, J., Major P.}, and {\sc Tusn\'{a}dy, G.} (1976).
An approximation of partial sums of independent rv's and the sample df. I.
{\em Z. Wahrsch. verw. Gebiete} {\bf 32}, 111--131.

\item{\sc Kosorok, M. R.} (1999). Two-sample quantile tests under
general conditions. {\em Biometrika} {\bf 86}, 909--921.

\item{\sc Kosorok, M. R.}, and {\sc Ma, S.} (2005). Comment on
``Semilinear high-dimensional model for normalization of microarray data:
a theoretical analysis and partial consistency'' by J. Fan, H. Peng,
T. Huang. {\em Journal of the American Statistical Association}, In press.

\item{\sc Massart, P.} (1990). The tight constant in the 
Dvoretsky-Kiefer-Wolfowitz inequality. {\em Annals of Probability}
{\bf 18}, 1269--1283.

\item{\sc Spang, R., Blanchette, C., Zuzan, H., Marks, J., Nevins, J.},
and {\sc West, M.} (2001). Prediction and uncertainty in the analysis
of gene expression profiles. In {\it Proceedings of the German Conference
on Bioinformatics GCB 2001}. Eds. E. Wingender, R. Hofestdt and I. Liebich,
Braunschweig, 102--111.

\item{\sc Storey, J. D., Taylor, J.E.}, and {\sc Siegmund, E.} (2004).
Strong control, conservative point estimation and simultaneous
conservative consistency of false discover rates: A unified approach.
{\em Journal of the Royal Statistical Society}, Series B, {\bf 66},
187--205.

\item{\sc van der Laan, M. J.}, and {\sc Bryan, J.} (2001). Gene
expression analysis with the parametric bootstrap. {\em Biostatistics}
{\bf 2}, 445--461.

\item{\sc van der Vaart, A. W.}, and {\sc Wellner, J. A.} (1996).
{\em Weak Convergence and Empirical Processes: With Applications to
Statistics}. Springer, New York.

\item{\sc West, M.} (2003) Bayesian factor regression models in
the "large p, small n" paradigm. In: {\it Bayesian Statistics} {\bf 7},
Eds. J. M. Bernardo, M. J. Bayarri, A. P. Dawid, J. O. Berger,
D. Heckerman, A. F. M. Smith and M. West, 733--742. 
Oxford University Press, Oxford.

\item{\sc West, M., Blanchette, C., Dressman, H., Huang, E., Ishida, S.,
Spang, R., Zuzan, H., Olson, J. A. Jr., Marks, J. R.}, and {\sc Nevins, 
J. R.} (2001). Predicting the clinical status of human breast cancer by
using gene expression profiles. {\em Proceedings of the National 
Academy of Sciences} {\bf 98}, 11462--11467.

\item{\sc Yang, Y.H., Dudoit, S., Luu, P., and Speed, T.P.} (2001)
Normalization for cDNA Microarray Data.
{\it Microarrays: Optical Technologies and Informatics, Vol. 4266 of
Proceedings of SPIE}, 141--152.

\end{list}

\vspace{0.25in}

\sc\begin{singlespace}
\noindent\begin{tabular}{ll}
M. R. Kosorok & S. Ma \\
Departments of Statistics and &Department of Biostatistics\\
Biostastistics/Medical Informatics & University of Washington\\
University of Wisconsin-Madison& Building 29, Suite 310 \\
1300 University Avenue & 6200 NE 74th St.\\
Madison, Wisconsin 53706&Seattle, WA 98115\\  
E-mail: {\rm kosorok{@}biostat.wisc.edu}&
E-mail: {\rm shuangge@u.washington.edu}
\end{tabular}
\rm

\newpage
\noindent
Table~1. One sample simulation study results for the mean, median
and signed rank statistics under models~1, 2 and 3. Tot.: 
total count identified 
using FDR. Pos.: number of true positives identified using FDR. EFDR:
empirical FDR. 
\begin{center}
\begin{tabular}{ccccccc}\hline
   & \multicolumn{2}{c}{Mean} & \multicolumn{2}{c}{Median} &
\multicolumn{2}{c}{Signed rank} \\
Model & Tot.(Pos.) & EFDR & Tot.(Pos.) & EFDR & Tot.(Pos.) & EFDR \\  \hline
\multicolumn{7}{c}{Sample size = 20}\\
 1 & 64.7(33.9) & 0.47& 31.8(25.4) & 0.19& 15.5(15.5) & 0.01\\
 2 & 64.4(33.9) & 0.47& 31.6(25.3) & 0.19& 15.3(15.2) & 0.01\\
 3 & 64.0(33.9) & 0.46& 31.1(25.0) & 0.19& 15.2(15.1) & 0.01\\
\hline
\multicolumn{7}{c}{Sample size = 50}\\
 1 & 54.2(37.8) & 0.30& 38.7(32.9) & 0.15& 34.5(34.0) & 0.01\\
 2 & 53.7(37.4) & 0.29& 38.5(32.7) & 0.14& 34.2(33.8) & 0.01\\
 3 & 52.3(37.5) & 0.27& 38.2(32.5) & 0.14& 34.4(33.9) & 0.01\\
\hline
\end{tabular}
\end{center}

\clearpage
\noindent
Table~2. Two sample simulation study results for mean, median,
Wilcoxon and Kolmogorov-Smirnov (KS) statistics under model~4.
Tot.: total count identified using FDR.
Pos: number of true positives identified using FDR. 
EFDR: empirical FDR. 

\begin{center}
\begin{tabular}{cccccccc}\hline
 \multicolumn{2}{c}{Mean} & \multicolumn{2}{c}{Median} &
\multicolumn{2}{c}{Wilcoxon} & \multicolumn{2}{c}{KS} \\
Tot.(Pos.) & EFDR & Tot.(Pos.) & EFDR & Tot.(Pos.) & 
EFDR &Tot.(Pos.) & EFDR \\  \hline
\multicolumn{8}{c}{$n_1=n_2=10$}\\
47.3(21.5)& 0.54&  8.4(6.7) & 0.18&   
14.4(13.1)& 0.08&  2.6(2.4) & 0.08 \\
\multicolumn{8}{c}{$n_1=n_2=30$}\\
40.9(28.9)& 0.28& 21.2(19.8)& 0.06&
32.0(26.6)& 0.16& 23.7(21.5)& 0.09\\
\multicolumn{8}{c}{$n_1=n_2=60$}\\
43.4(33.3)& 0.22& 29.7(25.4)& 0.14 &
39.4(32.4)& 0.17& 32.1(28.0)& 0.12\\
\hline
\end{tabular}
\end{center}  

\clearpage  

\end{singlespace}

\end{document}